\pdfoutput=1

\documentclass[aos]{imsart}

\RequirePackage{amsthm,amsmath,amsfonts,amssymb,amscd,bbm,enumitem,hyperref}
\usepackage[mathscr]{euscript}
\RequirePackage[numbers]{natbib}
\RequirePackage{graphicx}
\RequirePackage{stackengine,scalerel}
\RequirePackage{etoolbox}
\RequirePackage{anyfontsize}
\RequirePackage{thmtools}

\startlocaldefs
\theoremstyle{plain}

\theoremstyle{remark}
\newtheorem{defn}{Definition}

\newtheorem{thm}{Theorem}
\newtheorem{coro}{Corollary}
\theoremstyle{definition}             
\newtheorem{remark}{Remark}

\newtheorem{example}{Example}

\makeatother

\DeclareFontShape{T1}{ptm}{m}{scit}{<-> ssub * ptm/m/it}{}

\DeclareMathOperator*{\argmin}{argmin}

\DeclareMathOperator*{\polylog}{polylog}
\DeclareMathOperator*{\supp}{\operatorname{Supp}}
\DeclareMathOperator*{\lebmeas}{\mathfrak{m}}

\newcommand{\Indicator}{\mathbbm{1}}
\renewcommand{\P}{\mathbbm{P}}
\newcommand{\E}{\mathbbm{E}}

\newcommand{\cov}{\mathbbm{C}\mathrm{ov}}

\newcommand{\bnu}{\boldsymbol{\nu}}

\newcommand{\A}{\mathcal{A}}
\newcommand{\C}{\mathcal{C}}
\newcommand{\calK}{\mathcal{K}}
\newcommand{\calP}{\mathcal{P}}
\newcommand{\Q}{\mathcal{Q}}

\newcommand{\W}{\mathcal{W}}
\newcommand{\X}{\mathcal{X}}
\newcommand{\Y}{\mathcal{Y}}

\renewcommand{\S}{\mathcal{S}}

\newcommand{\G}{\mathscr{G}}
\newcommand{\scrH}{\mathscr{H}}
\newcommand{\R}{\mathscr{R}}
\newcommand{\F}{\mathscr{F}}

\newcommand{\reals}{\mathbb{R}}

\newcommand{\norm}[1]{\lVert #1 \rVert}
\newcommand{\infnorm}[1]{\lVert #1 \rVert_{\infty}}

\newcommand{\pOrder}{\mathfrak{p}}

\newcommand{\bH}{\mathbf{H}}
\newcommand{\bQ}{\mathbf{Q}}
\newcommand{\bS}{\mathbf{S}}
\newcommand{\bT}{\mathbf{T}}

\newcommand{\bV}{\mathbf{V}}

\newcommand{\be}{\mathbf{e}}
\newcommand{\bp}{\mathbf{p}}
\newcommand{\bu}{\mathbf{u}}
\newcommand{\bv}{\mathbf{v}}
\newcommand{\bw}{\mathbf{w}}
\newcommand{\bx}{\mathbf{x}}

\newcommand{\bz}{\mathbf{z}}

\newcommand{\bbeta}{\boldsymbol{\beta}}
\newcommand{\bgamma}{\boldsymbol{\gamma}}
\newcommand{\bSigma}{\boldsymbol{\Sigma}}

\newcommand{\ttv}{\mathtt{v}}
\newcommand{\ttK}{\mathtt{K}}
\newcommand{\ttk}{\mathtt{k}}

\newcommand{\ttc}{\mathtt{c}}
\newcommand{\ttd}{\mathtt{d}}
\newcommand{\ttL}{\mathtt{L}}
\newcommand{\ttM}{\mathtt{M}}
\newcommand{\ttN}{\mathtt{N}}
\newcommand{\ttE}{\mathtt{E}}
\newcommand{\ttS}{\mathtt{S}}
\newcommand{\ttTV}{\mathtt{TV}}

\newcommand{\bbQ}{\mathbb{Q}}

\newcommand{\proj}{\mathtt{\Pi}_{0}}

\newcommand{\projreg}{\mathtt{\Pi}_2}

\renewcommand{\d}{\mathfrak{d}}

\endlocaldefs

\begin{document}

\begin{frontmatter}
  \title{Strong Approximations for Empirical Processes\\ Indexed by Lipschitz Functions}
  \runtitle{Strong Approximations for Empirical Processes}

  \begin{aug}
    \author[A]{\fnms{Matias D.}~\snm{Cattaneo}\ead[label=e1]%
    {cattaneo@princeton.edu}},
    \and
    \author[A]{\fnms{Ruiqi (Rae)}~\snm{Yu}\ead[label=e2]%
    {rae.yu@princeton.edu}}
    \address[A]{Department of Operations Research and Financial Engineering,
      Princeton University, \\
    \printead[presep={}]{e1,e2}}
  \end{aug}

  \begin{abstract}
    This paper presents new uniform Gaussian strong approximations for empirical processes indexed by classes of functions based on $d$-variate random vectors ($d\geq1$). First, a uniform Gaussian strong approximation is established for general empirical processes indexed by possibly Lipschitz functions, improving on previous results in the literature. In the setting considered by \cite{Rio_1994_PTRF}, and if the function class is Lipschitzian, our result improves the approximation rate $n^{-1/(2d)}$ to $n^{-1/\max\{d,2\}}$, up to a $\polylog(n)$ term, where $n$ denotes the sample size. Remarkably, we establish a valid uniform Gaussian strong approximation at the rate $n^{-1/2}\log n$ for $d=2$, which was previously known to be valid only for univariate ($d=1$) empirical processes via the celebrated Hungarian construction \citep{komlos1975approximation}. Second, a uniform Gaussian strong approximation is established for multiplicative separable empirical processes indexed by possibly Lipschitz functions, which addresses some outstanding problems in the literature \citep[Section 3]{chernozhukov2014gaussian}. Finally, two other uniform Gaussian strong approximation results are presented when the function class is a sequence of Haar basis based on quasi-uniform partitions. Applications to nonparametric density and regression estimation are discussed.
  \end{abstract}



\end{frontmatter}


\section{Introduction}\label{section:Introduction}

Let $\bx_i\in \X \subseteq \mathbb{R}^d$, $i=1,\dots,n$, be independent and identical distributed (i.i.d.) random vectors supported on a background probability space $(\Omega,\mathcal{F},\P)$. The classical empirical process is
\begin{equation}\label{eq:X-process}
	X_n(h) = \frac{1}{\sqrt{n}} \sum_{i=1}^n \big( h(\bx_i) - \E[h(\bx_i)] \big), \qquad h \in \scrH,
\end{equation}
where $\scrH$ is a possibly $n$-varying class of functions. Following the empirical process literature, and assuming $\scrH$ is ``nice'', the stochastic process $(X_n(h):h \in \scrH)$ is said to be Donsker if it converges in law as $n\to\infty$ to a Gaussian process in $\ell^\infty(\scrH)$, the space of uniformly bounded real functions on $\scrH$. This weak convergence result is typically denoted by
\begin{equation}\label{eq:weak convergence}
	X_n \rightsquigarrow Z, \qquad \text{in } \ell^\infty(\scrH),
\end{equation}
where $(Z(h):h \in \scrH)$ is a mean-zero Gaussian process with covariance $\E[Z(h_1)Z(h_2)]=\E[h_1(\bx_i)h_2(\bx_i)]-\E[h_1(\bx_i)]\E[h_2(\bx_i)]$ for all $h_1,h_2 \in \scrH$ when $\scrH$ is not $n$-varying, or its limit as $n\to\infty$ otherwise. See \cite{wellner2013weak} and \cite{Gine-Nickl_2016_Book} for textbook overviews.

A more challenging endeavour is to construct a uniform Gaussian strong approximation for the empirical process $X_n$. That is, if the background probability space is ``rich'' enough, or is otherwise properly enlarged, the goal is to construct a sequence of mean-zero Gaussian processes $(Z_n(h):h \in \scrH)$ with the same covariance structure as $X_n$ (i.e., $\E[X_n(h_1)X_n(h_2)]=\E[Z_n(h_1)Z_n(h_2)]$ for all $h_1,h_2 \in \scrH$) such that
\begin{equation}\label{eq:SA}
	\|X_n - Z_n \|_\scrH = \sup_{h\in\scrH} \big| X_n(h) - Z_n(h) \big| = O(\varrho_n),
	\qquad \text{almost surely (a.s.)},
\end{equation}
for a non-random sequence $\varrho_n\to0$ as $n\to\infty$. Such a refined approximation result is useful in a variety of contexts. For example, it gives a distributional approximation for non-Donsker empirical processes, for which \eqref{eq:weak convergence} does not hold, and it also offers a precise quantification of the quality of the distributional approximation when \eqref{eq:weak convergence} holds. In addition, \eqref{eq:SA} is typically established using non-asymptotic probability concentration inequalities, which can be used to construct statistical inference procedures requiring uniformity over $\scrH$ and/or the class of underlying data generating processes. Furthermore, because the Gaussian process $Z_n$ is ``pre-asymptotic'', it can offer a better finite sample approximation to the sampling distribution of $X_n$ than the large sample approximation based on the limiting Gaussian process $Z$ in \eqref{eq:weak convergence}.

There is a large literature on strong approximations for empirical processes, offering different levels of tightness for the bound $\varrho_n$ in \eqref{eq:SA}. In particular, the univariate case ($d=1$) is mostly settled. A major breakthrough was accomplished by \cite[KMT hereafter]{komlos1975approximation}, who introduced the celebrated Hungarian construction to prove the optimal result $\varrho_n = n^{-1/2}\log n$ for the special case of the uniform empirical distribution process: $\bx_i\thicksim\mathsf{Uniform}(\X)$, $\X=[0,1]$, and $\scrH=\{\Indicator(\cdot\leq x):x\in[0,1]\}$, where $\Indicator(\cdot)$ denotes the indicator function. See \cite{Bretagnolle-Massart_1989_AOP} and \cite{Mason-VanZwet_2011_BookCh} for more technical discussions on the Hungarian construction, and \cite{csorgo1981strong}, \cite{Lindvall_1992_Book} and \cite{pollard2002user} for textbook overviews. The KMT result was later extended by \cite{Gine-Koltchinskii-Sakhanenko_2004_PTRF} and \cite{Gine-Nickl_2010_AoS} to univariate empirical processes indexed by functions with uniformly bounded total variation: for $\bx_i\thicksim\P_X$ supported on $\X=\mathbb{R}$ and continuously distributed, the authors obtained
\begin{equation}\label{eq:KMT-rate}
	\varrho_n = n^{-1/2} \log n,
\end{equation}
in \eqref{eq:SA}, with $\scrH$ satisfying a bounded variation condition. More recently, \cite[Lemma SA26]{Cattaneo-Feng-Underwood_2024_JASA} gave a self-contained proof of a slightly generalized KMT result allowing for a larger class of distributions $\P_X$. See Remark \ref{remark: comparison of tv} for details. As a statistical application, the authors considered univariate kernel density estimation \citep{Wand-Jones_1995_Book}, with bandwidth $b\to0$ as $n\to\infty$, and demonstrated that the optimal univariate KMT strong approximation rate $(nb)^{-1/2} \log n$ is achievable, where $nb$ is the effective sample size.

Establishing strong approximations for general empirical processes with $d\geq2$ is more difficult, since the KMT approach does not easily generalize to multivariate data. Foundational results include \cite{massart1989strong}, \cite{koltchinskii1994komlos}, and \cite{Rio_1994_PTRF}. In particular, assuming the function class $\scrH$ is uniformly bounded, has bounded total variation, and satisfies a VC-type condition, among other regularity conditions discussed precisely in the upcoming sections, \cite{Rio_1994_PTRF} obtained
\begin{equation}\label{eq:Rio-rate}
	\varrho_n = n^{-1/(2d)} \sqrt{\log n}, \qquad d\geq 2,
\end{equation}
in \eqref{eq:SA}. This result is tight under the conditions imposed \citep{beck1985lower}, and demonstrates an unfortunate dimension penalty in the convergence rate of the $d$-variate uniform Gaussian strong approximation. As a statistical application, the author also considered the kernel density estimator with bandwidth $b\to0$ as $n\to\infty$, and established \eqref{eq:SA} with
\[\varrho_n = (nb^d)^{-1/(2d)} \sqrt{\log n}, \qquad d\geq 2,\]
where $nb^d$ is the effective sample size.

While \cite{Rio_1994_PTRF}'s KMT strong approximation result is unimprovable under the conditions he imposed, it has two limitations:

\begin{enumerate}[leftmargin=*]
	\item The class of functions $\scrH$ may be too large, and further restrictions can open the door for improvements. For example, in his application to kernel density estimation, \cite[Section 4]{Rio_1994_PTRF} assumed that the class $\scrH$ is Lipschitzian to verify the sufficient conditions of his strong approximation theorem, but his theorem did not exploit the Lipschitz property in itself. (The Lipschitzian assumption is essentially without loss of generality in the kernel density estimation application.) It is an open question whether the optimal univariate KMT strong approximation rate \eqref{eq:KMT-rate} is achievable when $d\geq2$, under additional restrictions on $\scrH$.
	
	\item As discussed by \cite[Section 3]{chernozhukov2014gaussian}, applying \cite{Rio_1994_PTRF}'s strong approximation result directly to nonparametric local smoothing regression, a ``local empirical process'' in their terminology, leads to an even more suboptimal strong approximation rate in \eqref{eq:SA}. For example, in the case of kernel regression estimation with $d$-dimensional covariates, \cite{Rio_1994_PTRF}'s strong approximation would treat all $d+1$ variables (covariates and outcome) symmetrically, and thus it will give a strong approximation rate in \eqref{eq:SA} of the form
	\begin{align}\label{eq: local empirical process -- Rate Rio}
		\varrho_n = (n b^{d+1})^{-1/(2d+2)} \sqrt{\log n}, \qquad d\geq 1,
	\end{align}
	where $b\to0$ as $n\to\infty$, and under standard regular conditions. The main takeaway is that the resulting effective sample size is now $nb^{d+1}$ when in reality it should be $nb^{d}$, since only the $d$-dimensional covariates are smoothed out for estimation of the conditional expectation. It is this unfortunate fact that prompted \cite{chernozhukov2014gaussian} to develop strong approximation methods that target the scalar suprema of the stochastic process, $\sup_{h\in\scrH}|X_n(h)|$, instead of the stochastic process itself, $(X_n(h):h\in\scrH)$, as a way to circumvent the suboptimal strong approximation rates that would emerge from deploying directly \cite{Rio_1994_PTRF}'s result.
\end{enumerate}

This paper presents new uniform Gaussian strong approximation results for empirical processes that address the two aforementioned limitations. Section \ref{section:General Empirical Processes} studies the general empirical process \eqref{eq:X-process}, and establishes a uniform Gaussian strong approximation explicitly allowing for the possibility that $\scrH$ is Lipschitzian (Theorem \ref{thm: X-Process -- Main Theorem}). This result not only encompasses, but also generalizes previous results in the literature by allowing for $d\geq1$ under more generic entropy conditions and weaker conditions on the underlying data generating process. For comparison, if we impose the regularity conditions in \cite{Rio_1994_PTRF} and also assume $\scrH$ is Lipschitzian, then our result (Corollary \ref{coro: X-process -- vc lip}) verifies \eqref{eq:SA} with
\[\varrho_n = n^{-1/d} \sqrt{\log n} + n^{-1/2} \log n, \qquad d\geq 1,\]
thereby improving \eqref{eq:Rio-rate}, in addition to matching \eqref{eq:KMT-rate} when $d=1$; see Remark \ref{remark: comparison of tv} for details. Remarkably, we demonstrate that the optimal univariate KMT strong approximation rate $n^{-1/2} \log n$ is achievable when $d=2$, in addition to achieving the better approximation rate $n^{-1/d} \sqrt{\log n}$ when $d\geq3$. Applying our result to the kernel density estimation example, we obtain the improved strong approximation rate $(nb^d)^{-1/d} \sqrt{\log n} + (nb^d)^{-1/2} \log n$, $d\geq1$, under the same conditions imposed in prior literature. We thus show that the optimal univariate KMT uniform Gaussian strong approximation holds in \eqref{eq:SA} for bivariate kernel density estimation. Theorem \ref{thm: X-Process -- Main Theorem} also allows for other entropy notions for $\scrH$ beyond the classical VC-type condition, and delivers improvements over \cite{koltchinskii1994komlos}. See Remark \ref{remark: C^q -- pol entropy} for details. Section \ref{section:General Empirical Processes} discusses how our improvements are achieved, and outstanding technical roadblocks.

Section \ref{section: Residual Empirical Process} is motivated by the second aforementioned limitation in prior uniform Gaussian strong approximation results, and thus studies the \textit{residual-based empirical process}:
\begin{equation}\label{eq:R-process}
	R_n(g,r) = \frac{1}{\sqrt{n}} \sum_{i=1}^n \big( g(\bx_i) r(y_i) - \E[g(\bx_i) r(y_i)|\bx_i] \big),
	\qquad (g,r)\in \G\times\R,
\end{equation}
for $\bz_i = (\bx_i, y_i)$, $i=1,\dots, n$, a random sample now also including an outcome variable $y_i\in\reals$. Our terminology reflects the fact that $g(\bx_i) r(y_i) - \E[g(\bx_i) r(y_i)|\bx_i] = g(\bx_i) \epsilon_i(r)$ with $\epsilon_i(r) = r(y_i) - \E[r(y_i)|\bx_i]$, which can be interpreted as a residual in nonparametric local smoothing regression settings. In statistical applications, $g(\cdot)$ is typically an $n$-varying local smoother based on kernel, series, or nearest-neighbor methods, while $r(\cdot)$ is some transformation such as $r(y)=y$ for conditional mean or $r(y)=\Indicator(y\leq \cdot)$ for conditional distribution estimation. \cite[Section 3.1]{chernozhukov2014gaussian} call these special cases of $R_n$ a local empirical process.

The residual-based empirical process $(R_n(g,r): (g,r) \in \G \times \R)$ may be viewed as a general empirical process \eqref{eq:X-process} based on the sample $(\bz_i: 1 \leq i \leq n)$, and thus available strong approximation results can be applied directly, including \cite{koltchinskii1994komlos}, \cite{Rio_1994_PTRF}, and our new Theorem \ref{thm: X-Process -- Main Theorem}. However, those off-the-shelf results require stringent assumptions and can deliver suboptimal approximation rates. First, available results require $\bz_i$ to admit a bounded and positive Lebesgue density on $[0,1]^{d+1}$, possibly after some specific transformation, thereby imposing strong restrictions on the marginal distribution of $y_i$. Second, available results can lead to the incorrect effective sample size for the strong approximation rate. For example, for a local empirical process where $g(\cdot)$ denotes $n$-varying local smoothing weights based on a kernel function with bandwidth $b\to0$ as $n\to\infty$, and $r(y)=y$, \cite{Rio_1994_PTRF} gives the approximation rate \eqref{eq: local empirical process -- Rate Rio}, and our refined Theorem \ref{thm: X-Process -- Main Theorem} for general empirical processes indexed by Lipschitz functions gives a uniform Gaussian strong approximation rate
\begin{align}\label{eq: local empirical process -- Rate Thm 1}
	\varrho_n = (n b^{d+1})^{-1/(d+1)}\sqrt{\log n} + (n b^d)^{-1/2}\log n,
\end{align}
where the effective sample size is still $n b^{d+1}$. This is suboptimal because $n b^d$ is the (pointwise) effective sample size for the kernel regression estimator.

A key observation underlying the potential suboptimality of strong approximation results for local regression empirical processes is that all components of $\bz_i = (\bx_i, y_i)$ are treated symmetrically. Thus, Section \ref{section: Residual Empirical Process} presents a novel uniform Gaussian strong approximation for the residual-based empirical process (Theorem \ref{thm: R-Process -- Main Theorem}), which explicitly exploits the multiplicative separability of $\scrH = \G \times \R$ and the possibly Lipschitz continuity of the function class, while also removing stringent assumptions imposed on the underlying data generating process. When applied to the local kernel regression empirical processes, our best result gives a uniform Gaussian strong approximation rate
\begin{align}\label{eq: R-process -- SA2}
	\varrho_n = (n b^d)^{-1/(d+2)} \sqrt{\log n} + (n b^d)^{-1/2} \log n,
\end{align}
thereby improving over both \cite{Rio_1994_PTRF} leading to \eqref{eq:Rio-rate}, and Theorem \ref{thm: X-Process -- Main Theorem} leading to \eqref{eq: local empirical process -- Rate Thm 1}. The correct effective sample size $n b^d$ is achieved, under weaker regularity conditions. As a statistical application, Section \ref{section: Local Polynomial Regression} leverages Theorem \ref{thm: R-Process -- Main Theorem} to establish the best known uniform Gaussian strong approximation result for local polynomial regression estimators \citep{Fan-Gijbels_1996_Book}. 

Following \cite{Rio_1994_PTRF}, the proof of Theorem \ref{thm: X-Process -- Main Theorem} in Section \ref{section:General Empirical Processes} first approximates in mean square the class of functions $\scrH$ using a Haar basis over carefully constructed disjoint \textit{dyadic} cells, and then applies the celebrated Tusn\'ady's Lemma \cite[][Chapter 10, for a textbook introduction]{pollard2002user} to construct a strong approximation. It thus requires balancing two approximation errors: a projection error (``bias'') emerging from the mean square projection based on a Haar basis, and a coupling error (``variance'') emerging from the coupling construction for the projected process. A key observation in our paper is that both errors can be improved by explicitly exploiting a Lipschitz assumption on $\scrH$. However, it appears that to achieve the univariate KMT uniform Gaussian strong approximation for the general empirical process \eqref{eq:X-process} with $d\geq3$, a mean square projection based on a higher-order function class would be needed to improve the projection error, but no coupling methods available in the literature for the resulting projected process. The proof of Theorem \ref{thm: R-Process -- Main Theorem} in Section \ref{section: Residual Empirical Process} employs a similar projection and coupling decomposition approach, but treats $\G$ and $\R$ separately in order to leverage the multiplicative separability of the residual process $(R_n(g,r): (g,r) \in \G \times \R)$. In particular, the proof designs cells for projection and coupling approximation that are asymmetric in the direction of $\bx_i$ and $y_i$ components to obtain the uniform Gaussian strong approximation. This distinct proof strategy relaxes some underlying assumptions (most notably, on the distribution of $y_i$), and delivers a better strong approximation rate for some local empirical processes than what would be obtained by directly applying Theorem \ref{thm: X-Process -- Main Theorem}.   

In general, however, neither Theorem \ref{thm: X-Process -- Main Theorem} nor Theorem \ref{thm: R-Process -- Main Theorem} dominates each other, nor their underlying assumptions imply each other, and therefore both are of interest, depending on the statistical problem under consideration. Their proofs employ different strategies (most notably, in terms of the dyadic cells expansion used) to leverage the specific structure of $X_n$ and $R_n$. It is an open question whether the uniform Gaussian strong approximation rates obtained from Theorems \ref{thm: X-Process -- Main Theorem} and \ref{thm: R-Process -- Main Theorem} are optimal under the assumptions imposed.

As a way to circumvent the technical limitations underlying the proof strategies of Theorem \ref{thm: X-Process -- Main Theorem} and Theorem \ref{thm: R-Process -- Main Theorem}, Section \ref{sec: Quasi-Uniform Haar Basis} presents two other uniform Gaussian strong approximation results when $\scrH$ is spanned by a possibly increasing sequence of finite Haar functions on \textit{quasi-uniform} partitions, for the general empirical process (Theorem \ref{thm: X-Process -- Approximate Dyadic Thm}) and for the residual-based empirical process (Theorem \ref{thm: R-Process -- Approximate Dyadic Thm}). These theorems shut down the projection error, and also rely on a generalized Tusn\'ady's Lemma established in this paper, to establish valid couplings over more general partitioning schemes and under weaker regularity conditions. In this specialized setting, we demonstrate that a uniform Gaussian strong approximation at the optimal univariate KMT rate based on the corresponding effective sample size is possible for all $d\geq1$, up to a $\polylog(n)$ term, where $\polylog(n) = \log^a(n)$ for some $a>0$, and possibly an additional ``bias'' term induced exclusively by the cardinality of $\R$. As statistical applications, we establish uniform Gaussian strong approximations for the classical histogram density estimator, and for Haar partitioning-based regression estimators such as those arising in the context of certain regression tree and related nonparametric methods \citep{breiman1984,Huang_2003_AOS,Cattaneo-Farrell-Feng_2020_AOS}.

The supplemental appendix \cite{Cattaneo-Yu_2024_AOS--SA} contains all technical proofs, additional theoretical results of independent interest, and other omitted details.

\subsection{Related Literature}

This paper contributes to the literature on strong approximations for empirical processes, and their applications to uniform inference for nonparametric smoothing methods. For introductions and overviews, see \cite{csorgo1981strong}, \cite{Lindvall_1992_Book}, \cite{Einmahl-Mason_1998_BookCh}, \cite{berthet2006revisiting}, \cite{Mason-Zhou_2012_PS}, \cite{Gine-Nickl_2016_Book}, \cite{pollard2002user}, \cite{Zaitsev_2013_RMSurveys}, and references therein. See also \cite[Section 3]{chernozhukov2014gaussian} for discussion and further references concerning local empirical processes and their role in nonparametric curve estimation.

The celebrated KMT construction \cite{komlos1975approximation}, Yurinskii's coupling \cite{yurinskii1978error}, and Zaitsev's coupling \cite{zauitsev1987estimates} are three well-known approaches that can be used to establish a uniform Gaussian strong approximation for empirical processes. Among them, the KMT approach often offers the tightest approximation rates when applicable, and is the focus of our paper: closely related literature includes \cite{massart1989strong}, \cite{koltchinskii1994komlos}, \cite{Rio_1994_PTRF}, \cite{Gine-Koltchinskii-Sakhanenko_2004_PTRF}, and \cite{Gine-Nickl_2010_AoS}, among others. As summarized in the introduction, our first main result (Theorem \ref{thm: X-Process -- Main Theorem}) encompasses and improves on prior results in that literature. Furthermore, Theorems \ref{thm: R-Process -- Main Theorem}, \ref{thm: X-Process -- Approximate Dyadic Thm}, and \ref{thm: R-Process -- Approximate Dyadic Thm} offer new results for more specific settings of interest in statistics, in particular addressing some outstanding problems in the literature \citep[Section 3]{chernozhukov2014gaussian}. We provide detailed comparisons to the prior literature in the upcoming sections.

We do not discuss the other coupling approaches because they deliver slower strong approximation rates under the assumptions imposed in this paper: for example, see \cite{cattaneo2024yurinskii} for results based on Yurinskii's coupling, and \cite{settati2009gaussian} for results based on Zaitsev's coupling. Finally, employing a different approach, \cite{dedecker2014strong} obtain a uniform Gaussian strong approximation for the multivariate empirical process indexed by half plane indicators with a dimension-independent approximation rate, up to $\polylog(n)$ terms.

\section{Notation}\label{section:Setup}

We employ standard notations from the empirical process literature, suitably modified and specialized to improve exposition. See, for example, \cite{Ambrosio-Fusco-Pallara_2000_book}, \cite{wellner2013weak} and \cite{Gine-Nickl_2016_Book} for background definitions and more details.

The $q$-dimensional Gaussian distribution with mean $\boldsymbol{\mu} \in \reals^q$ and symmetric positive semidefinite covariance matrix $\boldsymbol{\Sigma} \in \reals^{q \times q}$ is denoted by $\mathsf{Normal}(\boldsymbol{\mu},\boldsymbol{\Sigma})$. The binomial distribution with parameter $n \in \mathbb{N}$ and $p \in [0,1]$ is denoted by $\mathsf{Bin}(n,p)$. $|\A|$ denotes the cardinality of the set $\A$. For a vector $\mathbf{a} \in \reals^{q}$, $\norm{\mathbf{a}}$ denotes the Euclidean norm and $\infnorm{\mathbf{a}}$ denotes the maximum norm of $\mathbf{a}$. For a matrix $\mathbf{A} \in \reals^{q \times q}$, $\sigma_1(\mathbf{A}) \geq \sigma_2(\mathbf{A}) \geq \cdots \geq \sigma_d(\mathbf{A}) \geq 0$ denote the singular values of $\mathbf{A}$. For $1 \leq i_1 \leq j_2 \leq n$ and $1 \leq j_1 \leq j_2 \leq n$, $\mathbf{A}_{i_1:i_2,j_1:j_2}$ denotes the submatrix $(A_{ij})_{i_1 \leq i \leq i_2,j_1 \leq j \leq j_2}$ of $\mathbf{A}$, and $\mathbf{A}_{i_1,j_1:j_2}$, $\mathbf{A}_{i_1:i_2,j_1}$ are likewise defined. For sequences of real numbers, we write $a_n = o(b_n)$ if $\limsup_{n \to \infty} |\frac{a_n}{b_n}| = 0$, and write $a_n=O(b_n)$ if there exists some constant $C$ and $N > 0$ such that $n > N$ implies $|a_n| \leq C |b_n|$. For sequences of random variables, we write $a_n = o_{\P}(b_n)$ if $\limsup_{n \to \infty}\P[|\frac{a_n}{b_n}| \geq \varepsilon] = 0$ for all $\varepsilon > 0$, and write $a_n=O_\P(b_n)$ if $\limsup_{M \to \infty} \limsup_{n \to \infty} \P[|\frac{a_n}{b_n}| \geq M] = 0$.

Let $\mathcal{U}, \mathcal{V} \subseteq \reals^q$. We define $\mathcal{U} + \mathcal{V} = \{\bu + \bv: \bu \in \mathcal{U}, \bv \in \mathcal{V}\}$ and $\infnorm{\mathcal{U}} = \sup\{\infnorm{\bu_1 - \bu_2}: \bu_1,\bu_2 \in \mathcal{U}\}$, and $\mathcal{B}(\mathcal{U})$ denotes the Borel $\sigma$-algebra generated by $\mathcal{U}$ and $\mathcal{B}(\mathcal{U}) \otimes \mathcal{B}(\mathcal{V})$ denotes the product $\sigma$-algebra. Let $\mu$ be a measure on $(\mathcal{U},\mathcal{B}(\mathcal{U}))$, and $\phi: (\mathcal{V},\mathcal{B}(\mathcal{V})) \mapsto (\mathcal{U},\mathcal{B}(\mathcal{U}))$ be a measurable function. $\mu \circ \phi$ denotes the measure on $(\mathcal{V},\mathcal{B}(\mathcal{V}))$ such that $\mu \circ \phi(V) = \mu(\phi(V))$ for any $V \in \mathcal{B}(\mathcal{V})$. For $R \in \mathcal{B}(\mathcal{U})$, let $\mu|_R$ be the restriction of $\mu$ on $R$, that is, $\mu|_R(U) = \mu(U \cap R)$ for all $U \in \mathcal{B}(\mathcal{U})$. Two measures $\mu$ and $\nu$ on the measure space $(\mathcal{U},\mathcal{B}(\mathcal{U}))$ agree on $R\in \mathcal{B}(\mathcal{U})$ if $\mu|_R = \nu|_R$. The support of $\mu$ is $\supp(\mu)=\operatorname{closure}(\cup\{U \in \mathcal{B}(\mathcal{U}): \mu(U) \neq 0\})$. The Lebesgue measure is denoted by $\lebmeas$. Let $f$ be a real-valued function on the measure space $(\mathcal{U},\mathcal{B}(\mathcal{U}),\mu)$. Define the $L_p$ norms $\norm{f}_{\mu,p} =(\int |f|^p d \mu)^{1/p}$ for $1 \leq p < \infty$ and $\norm{f}_{\infty} = \sup_{\bx \in \mathcal{U}} |f(\bx)|$, and let $\supp(f) = \{\bu \in \mathcal{U}: f(\bu) > 0\}$ be the support of $f$. $L_p(\mu)$ is the class of all real-valued measurable functions $f$ on $(\mathcal{U}, \mathcal{B}(\mathcal{U}))$ such that $\|f\|_{\mu,p}< \infty$, for $1 \leq p < \infty$. The semi-metric $\d_\mu$ on $L_2(\mu)$ is defined by $\d_\mu(f, g) = (\|f - g\|_{\mu,2}^2 - (\int f \, d\mu - \int g \, d\mu)^2)^{1/2}$, for $f, g \in L_2(\mu)$. Whenever it exits, $\nabla f(\bx)$ denotes the Jacobian matrix of $f$ at $\bx$. If $\F$ and $\G$ are two sets of functions from measure space $(\mathcal{U},\mathcal{B}(\mathcal{U}))$ and $(\mathcal{V},\mathcal{B}(\mathcal{V}))$ to $\reals$, respectively, then $\F \times \G$ denotes the class of measurable functions $\{(f,g): f \in \F, g \in \G\}$ from $(\mathcal{U} \times \mathcal{V}, \mathcal{B}(\mathcal{U}) \otimes \mathcal{B}(\mathcal{V}))$ to $\reals$. For a measure $\mu$ on $(\mathcal{U} \times \mathcal{V},\mathcal{B}(\mathcal{U}) \otimes \mathcal{B}(\mathcal{V}))$, the semi-metric $\d_{\mu}$ on $\G \times \R$ is defined by $\d_{\mu}((g_1,r_1),(g_2,r_2)) = (\norm{g_1 r_1 - g_2 r_2}_{\mu,2}^2 - (\int g_1 r_1 d \mu - \int g_2 r_2 d \mu)^2)^{1/2}$. For a semi-metric space $(\S,d)$, the covering number $N(\S, d, \varepsilon)$ is the minimal number of balls $B_v(\varepsilon) =\{u: d(u,v) < \varepsilon\}$, $v\geq1$, needed to cover $\S$. 

\subsection{Main Definitions}\label{section:Main Definitions}

Let $\F$ be a class of measurable functions from a probability space $(\reals^q, \mathcal{B}(\reals^q), \P)$ to $\reals$. We introduce several definitions that capture properties of $\F$.

\begin{defn}\label{sa-defn: meas}
	$\F$ is pointwise measurable if it contains a countable subset $\G$ such that for any $f \in \F$, there exists a sequence $(g_m:m\geq1) \subseteq \G$ such that $\lim_{m \to \infty} g_m(\bu) = f(\bu)$ for all $\bu \in \reals^q$.
\end{defn}

\begin{defn}\label{defn: normalizing transformation}
	Let $\supp(\F) = \cup_{f \in \F}\supp(f)$. A probability measure $\bbQ_\F$ on $(\reals^q,\mathcal{B}(\reals^q))$ is a surrogate measure for $\P$ with respect to $\F$ if
	\begin{enumerate}[label=(\roman*)]
		\item $\bbQ_\F$ agrees with $\P$ on $\supp(\P) \cap \supp(\F)$.
		\item $\bbQ_\F(\supp(\F) \setminus \supp(\P)) = 0$.
	\end{enumerate}
	Let $\Q_\F=\supp(\bbQ_\F)$.
\end{defn}

\begin{defn}\label{sa-defn: pointwise tv}
	For $q=1$ and an interval $\mathcal{I}\subseteq\reals$, the pointwise total variation of $\F$ over $\mathcal{I}$ is 
	\begin{align*}
		\mathtt{pTV}_{\F,\mathcal{I}} = \sup_{f \in \F} \sup_{P\geq1}\sup_{\calP_P \in \mathcal{I}} \sum_{i = 1}^{P-1}|f(a_{i+1}) - f(a_i)|,
	\end{align*}
	where $\calP_P=\{(a_1,\dots,a_P):a_1 \leq \cdots \leq a_P\}$ denotes the collection of all partitions of $\mathcal{I}$.
\end{defn}

\begin{defn}\label{defn: tv}
	For a non-empty $\mathcal{C} \subseteq \reals^q$, the total variation of $\F$ over $\mathcal{C}$ is 
	\begin{align*}
		\ttTV_{\F, \mathcal{C}} = \inf_{\mathcal{U} \in \mathcal{O}(\mathcal{C})}\sup_{f \in \F} \sup_{\phi \in \mathscr{D}_{q}(\mathcal{U})} \int_{\reals^q} f(\bu)\operatorname{div}(\phi)(\bu) d \bu / \infnorm{\norm{\phi}_2},
	\end{align*}
	where $\mathcal{O}(\mathcal{C})$ denotes the collection of all open sets that contains $\mathcal{C}$, and $\mathscr{D}_{q}(\mathcal{U})$ denotes the space of infinitely differentiable functions from $\reals^q$ to $\reals^q$ with compact support contained in $\mathcal{U}$. 
\end{defn}

\begin{defn}\label{sa-defn: local tv}
	For a non-empty $\mathcal{C} \subseteq \reals^q$, the local total variation constant of $\F$ over $\mathcal{C}$, is a positive number $\ttK_{\F,\mathcal{C}}$ such that for any cube $\mathcal{D} \subseteq \reals^q$ with edges of length $\ell$ parallel to the coordinate axises, 
	\begin{align*}
		\ttTV_{\F, \mathcal{D} \cap \mathcal{C}}  \leq \ttK_{\F, \mathcal{C}} \ell^{d-1}.
	\end{align*}
\end{defn}

\begin{defn}\label{sa-defn: envelope}
	For a non-empty $\mathcal{C} \subseteq \reals^q$, the envelopes of $\F$ over $\mathcal{C}$ are
	\begin{align*}
		\ttM_{\F,\C} = \sup_{\bu \in \C }M_{\F,\C}(\bu),
		\qquad M_{\F,\C}(\bu) = \sup_{f \in \F}|f(\bu)|,
		\qquad \bu \in \C.
	\end{align*}
\end{defn}

\begin{defn}\label{sa-defn: Lipschitz}
	For a non-empty $\mathcal{C} \subseteq \reals^q$, the Lipschitz constant of $\F$ over $\C$ is
	\begin{align*}
		\ttL_{\F,\C} = \sup_{f \in \F}\sup_{\bu_1, \bu_2 \in \C} \frac{|f(\bu_1) - f(\bu_2)|}{\|\bu_1 - \bu_2\|_\infty}.
	\end{align*}
\end{defn}

\begin{defn}\label{sa-defn: L1}
	For a non-empty $\mathcal{C} \subseteq \reals^q$, the $L_1$ bound of $\F$ over $\C$ is
	\begin{align*}
		\mathtt{E}_{\F,\C} = \sup_{f \in \F} \int_{\C} |f| d\P.
	\end{align*}
\end{defn}

\begin{defn}\label{sa-defn: uniform covering number}
	For a non-empty $\mathcal{C} \subseteq \reals^q$, the uniform covering number of $\F$ with envelope $M_{\F,\C}$ over $\C$ is
	\begin{align*}
		\ttN_{\F,\C}(\delta,M_{\F,\C}) = \sup_{\mu} N(\F,\norm{\cdot}_{\mu,2},\delta \norm{M_{\F,\C}}_{\mu,2}),
		\qquad \delta \in (0, \infty),
	\end{align*}
	where the supremum is taken over all finite discrete measures on $(\C, \mathcal{B}(\C))$. We assume that $M_{\F,\C}(\bu)$ is finite for every $\bu \in \C$.
\end{defn}

\begin{defn}\label{sa-defn: uniform entropy integrals}
	For a non-empty $\mathcal{C} \subseteq \reals^q$, the uniform entropy integral of $\F$ with envelope $M_{\F,\C}$ over $\C$ is
	\begin{align*}
		J_\C(\delta, \F, M_{\F,\C}) = \int_0^{\delta} \sqrt{1 + \log \ttN_{\F,\C}(\varepsilon,M_{\F,\C})} d \varepsilon,
	\end{align*}
	where it is assumed that $M_{\F,\C}(\bu)$ is finite for every $\bu \in \C$.
\end{defn}

\begin{defn}\label{sa-defn: vc}
	For a non-empty $\mathcal{C} \subseteq \reals^q$, $\F$ is a VC-type class with envelope $M_{\F,\C}$ over $\C$ if (i) $M_{\F,\C}$ is measurable and $M_{\F,\C}(\bu)$ is finite for every $\bu \in \C$, and (ii) there exist $\ttc_{\F,\C}>0$ and $\mathtt{d}_{\F,\C}>0$ such that
	\begin{align*}
		\ttN_{\F,\C}(\varepsilon,M_{\F,\C}) \leq \ttc_{\F,\C} \varepsilon^{-\mathtt{d}_{\F,\C}},
		\qquad \varepsilon\in(0,1).
	\end{align*}
\end{defn}

\begin{defn}\label{sa-defn: ploynomial-entropy}
	For a non-empty $\mathcal{C} \subseteq \reals^q$, $\F$ is a polynomial-entropy class with envelope $M_{\F,\C}$ over $\C$ if (i) $M_{\F,\C}$ is measurable and $M_{\F,\C}(\bu)$ is finite for every $\bu \in \C$, and (ii) there exist $\mathtt{a}_{\F,\C}>0$ and $\mathtt{b}_{\F,\C}>0$ such that
	\begin{align*}
		\log \ttN_{\F,\C}(\varepsilon,M_{\F,\C}) \leq \mathtt{a}_{\F,\C} \varepsilon^{-\mathtt{b}_{\F,\C}},
		\qquad \varepsilon\in(0,1).
	\end{align*}
\end{defn}

If a surrogate measure $\bbQ_\F$ for $\P$ with respect to $\F$ has been assumed, and it is clear from the context, we drop the dependence on $\C = \Q_{\F}$ for all quantities in Definitions \ref{defn: tv}--\ref{sa-defn: ploynomial-entropy}. That is, to save notation, we set $\ttTV_{\F}=\ttTV_{\F,\Q_{\F}}$, $\ttK_{\F}=\ttK_{\F,\Q_{\F}}$, $\ttM_{\F}=\ttM_{\F,\Q_{\F}}$, $M_{\F}(\bu)=M_{\F,\Q_{\F}}(\bu)$, $\ttL_{\F}=\ttL_{\F,\Q_{\F}}$, and so on, whenever there is no confusion.

\section{General Empirical Process}\label{section:General Empirical Processes}

Let
\begin{equation*}
	\mathsf{m}_{n,d} =              
	\begin{cases}
		n^{-1/2}\sqrt{\log n} & \text{ if } d = 1\\
		n^{-1/(2d)}           & \text{ if } d \geq 2
	\end{cases} \qquad\text{and}\qquad
	\mathsf{l}_{n,d} =
	\begin{cases}
		1                     & \text{ if } d = 1\\
		n^{-1/2}\sqrt{\log n} & \text{ if } d = 2\\
		n^{-1/d}              & \text{ if } d \geq 3
	\end{cases},
\end{equation*}
and recall Section \ref{section:Main Definitions} and the notation conventions introduced there.

\begin{thm}\label{thm: X-Process -- Main Theorem}
	Suppose $(\bx_i: 1 \leq i \leq n)$ are i.i.d. random vectors taking values in $(\reals^d, \mathcal{B}(\reals^d))$ with common law $\P_X$ supported on $\X \subseteq \reals^d$, and the following conditions hold.
	\begin{enumerate}[label=(\roman*)]
		\item $\scrH$ is a real-valued pointwise measurable class of functions on $(\reals^d,\mathcal{B}(\reals^d),\P_X)$.
		\item There exists a surrogate measure $\bbQ_\scrH$ for $\P_X$ with respect to $\scrH$ such that $\bbQ_\scrH = \lebmeas \circ \phi_\scrH$, where the \textit{normalizing transformation} $\phi_{\scrH}: \Q_\scrH \mapsto [0,1]^d$ is a diffeomorphism.
		\item $\ttM_{\scrH} < \infty$ and $J(1, \scrH, \ttM_{\scrH}) < \infty$.
	\end{enumerate}
	
	Then, on a possibly enlarged probability space, there exists a sequence of mean-zero Gaussian processes $(Z^X_n(h):h\in\scrH)$ with almost sure continuous trajectories on $(\scrH, \d_{\P_X})$ such that:
	\begin{itemize}
		\item $\E[X_n(h_1) X_n(h_2)] = \E[Z^X_n(h_1) Z^X_n(h_2)]$ for all $h_1, h_2 \in \scrH$, and
		\item $\P\big[\norm{X_n - Z^X_n}_{\scrH} > C_1 \mathsf{S}_n(t)\big] \leq C_2 e^{-t}$ for all $t > 0$,
	\end{itemize}
	where $C_1$ and $C_2$ are universal constants, and
	\[\mathsf{S}_n(t) = \min_{\delta\in(0,1)}\{\mathsf{A}_n(t,\delta)+\mathsf{F}_n(t,\delta)\},\]
	where
	\begin{align*}
		\mathsf{A}_n(t,\delta)
		&= \min\big\{ \mathsf{m}_{n,d} \sqrt{\ttM_{\scrH}}, \mathsf{l}_{n,d} \sqrt{\ttc_2\ttL_{\scrH}} \big\} \sqrt{\ttc_1 \ttTV_{\scrH}} \sqrt{t + \log \ttN_\scrH(\delta, M_\scrH)}\\
		&\qquad + \sqrt{\frac{\ttM_{\scrH}}{n}}  \min\big\{ \sqrt{\log n} \sqrt{\ttM_{\scrH}} , \sqrt{ \ttc_3 \ttK_{\scrH} + \ttM_\scrH}\big\} (t + \log \ttN_\scrH(\delta,M_\scrH))
	\end{align*}
	\begin{align*}
		\ttc_1 = d \sup_{\bx \in \Q_{\scrH}} \prod_{j = 1}^{d-1} \sigma_j(\nabla \phi_{\scrH}(\bx)),
		\quad
		\ttc_2 = \sup_{\bx \in \Q_{\scrH}} \frac{1}{\sigma_{d}(\nabla \phi_{\scrH}(\bx))},
		\quad
		\ttc_3 = 2^{d-1} d^{d/2-1} \ttc_1 \ttc_2^{d-1},
	\end{align*}
	and
	\begin{align*}
		\mathsf{F}_n(t,\delta) = J(\delta, \scrH, \ttM_{\scrH}) \ttM_{\scrH} + \frac{\ttM_{\scrH} J^2(\delta, \scrH, \ttM_{\scrH})}{\delta^2 \sqrt{n}} + \delta \ttM_{\scrH} \sqrt{t} + \frac{\ttM_{\scrH}}{\sqrt{n}} t.
	\end{align*}
\end{thm}

This uniform Gaussian strong approximation theorem is given in full generality to accommodate different applications. Section \ref{section: normalizing transformation} discusses the role of the surrogate measure and normalizing transformation, and Section \ref{sec: Special Cases of Theorem 1} discusses leading special cases and compares our results to prior literature. The proof of Theorem \ref{thm: X-Process -- Main Theorem} is in \cite[Section SA-II]{Cattaneo-Yu_2024_AOS--SA}, but we briefly outline the general proof strategy here to highlight our improvements on prior literature and some open questions. The proof begins with the standard discretization (or meshing) decomposition:
\begin{align*}
	\norm{X_n - Z_n^X}_{\scrH} \leq \norm{X_n - X_n\circ\pi_{\scrH_\delta}}_{\scrH} + \norm{X_n - Z_n^X}_{\scrH_\delta} + \norm{Z_n^X\circ\pi_{\scrH_\delta}-Z_n^X}_{\scrH},
\end{align*}
where $\norm{X_n - Z_n^X}_{\scrH_\delta}$ captures the coupling between the empirical process and the Gaussian process on a $\delta$-net of $\scrH$, which is denoted by $\scrH_{\delta}$, while the terms $\norm{X_n - X_n\circ\pi_{\scrH_\delta}}_{\scrH}$ and $\norm{Z_n^X\circ\pi_{\scrH_\delta}-Z_n^X}_{\scrH}$ capture the fluctuations (or oscillations) relative to the meshing for each of the stochastic processes. The latter two errors are handled using standard empirical process results, which give the contribution $\mathsf{F}_n(t,\delta)$ emerging from Talagrand's inequality \citep[Theorem 3.3.9]{Gine-Nickl_2016_Book} combined with a standard maximal inequality \citep[Theorem 5.2]{chernozhukov2014gaussian}. 

Following \cite{Rio_1994_PTRF}, the coupling term $\norm{X_n - Z_n^X}_{\scrH_\delta}$ is further decomposed using a mean square projection onto a Haar function space:
\begin{align}\label{eq: X-Process -- Coupling}
	\norm{X_n - Z_n^X}_{\scrH_\delta}
	\leq \norm{X_n - \proj X_n}_{\scrH_{\delta}} + \norm{\proj X_n - \proj Z_n^X}_{\scrH_{\delta}} + \norm{\proj Z_n^X - Z_n^X}_{\scrH_\delta},
\end{align}
where $\proj X_n(h) = X_n \circ \proj h$ with $\proj$ denoting the $L_2$-projection onto piecewise constant functions on a carefully chosen partition of $\X$. We introduce a class of recursive \textit{quasi-dyadic} cells expansion of $\X$, which we employ to generalize prior results in the literature, including properties of the $L_2$-projection onto a Haar basis based on quasi-dyadic cells.

The term $\norm{\proj X_n - \proj Z_n^X}_{\scrH_{\delta}}$ in \eqref{eq: X-Process -- Coupling} represents the strong approximation error for the projected process over a recursive dyadic collection of cells partitioning $\X$. Handling this error boils down to the coupling of $\mathsf{Bin}(n,\frac{1}{2})$ with $\mathsf{Normal}(\frac{n}{2}, \frac{n}{4})$, due to the fact that the constant approximation within each recursive partitioning cell generates counts based on i.i.d. data. Building on the celebrated Tusn\'ady's Lemma, \cite[Theorem 2.1]{Rio_1994_PTRF} established a remarkable coupling result for bounded functions $L_2$-projected on a dyadic cells expansion of $\X$. We build on his powerful ideas, and establish an analogous result for the case of Lipschitz functions $L_2$-projected on dyadic cells expansion of $\X$, thereby obtaining a tighter coupling error. A limitation of these results is that they only apply to a dyadic cells expansion due to the specifics of Tusn\'ady's Lemma.

The terms $\norm{X_n - \proj X_n}_{\scrH_{\delta}}$ and $\norm{\proj Z_n^X - Z_n^X}_{\scrH_\delta}$ in \eqref{eq: X-Process -- Coupling} represent the errors of the mean square projection onto a Haar basis based on \textit{quasi-dyadic} cells expansion of $\X$. We handle this error using Bernstein inequality, while also taking into account explicitly the potential Lipschitz structure of the functions, and the more generic cell structure.

Balancing the coupling error and the two projection errors in \eqref{eq: X-Process -- Coupling} gives term $\mathsf{A}_n(t,\delta)$ in Theorem~\ref{thm: X-Process -- Main Theorem}. Section SA-II of \cite{Cattaneo-Yu_2024_AOS--SA} provides all technical details, and additional results that may be of independent interest.

\subsection{Surrogate Measure and Normalizing Transformation}\label{section: normalizing transformation}

Theorem \ref{thm: X-Process -- Main Theorem} assumes the existence of a surrogate measure $\bbQ_\scrH$, and a normalizing transformation $\phi_{\scrH}$, which together restrict $\P_X$ to be absolutely continuous with respect to $\lebmeas$ on $\X \cap \supp(\scrH)$, while incorporating features of the support of $\scrH$. We provide examples of $\bbQ_\scrH$ and $\phi_{\scrH}$, discuss primitive sufficient conditions, and bound the constants $\ttc_1$, $\ttc_2$, and $\ttc_3$ explicitly.

As a first simple example, suppose that $\bx_i\thicksim\mathsf{Uniform}(\X)$ with $\X=\times_{l = 1}^d[\mathsf{a}_l, \mathsf{b}_l]$, where $-\infty < \mathsf{a}_l < \mathsf{b}_l < \infty$, $l=1,2,\dots,d$. Setting $\bbQ_\scrH=\P_X$ and $\phi_{\scrH}(x_1, \cdots, x_d) = ((\mathsf{b}_1 - \mathsf{a}_1)^{-1}(x_1 - \mathsf{a}_1), \cdots, (\mathsf{b}_d - \mathsf{a}_d)^{-1}(x_d - \mathsf{a}_d))$ verifies assumption (ii) in Theorem \ref{thm: X-Process -- Main Theorem}. In this case, $\ttc_1 = d \max_{1 \leq l \leq d}|\mathsf{b}_l - \mathsf{a}_l| \prod_{l = 1}^d |\mathsf{b}_l - \mathsf{a}_l|^{-1}$, $\ttc_2 = \max_{1 \leq l \leq d}|\mathsf{b}_l - \mathsf{a}_l|$ and $\ttc_3 = 2^{d-1} d^{d/2} \max_{1 \leq l \leq d}|\mathsf{b}_l - \mathsf{a}_l|^d \prod_{l = 1}^d |\mathsf{b}_l - \mathsf{a}_l|^{-1}$.

When $\P_X$ is not the uniform distribution, or $\X$ is not isomorphic to the $d$-dimensional unit cube, a careful choice of $\bbQ_{\scrH}$ and $\phi_{\scrH}$ is needed. In many interesting cases, the \emph{Rosenblatt transformation} can be used to exhibit a valid normalizing transformation, together with an appropriate choice of $\bbQ_{\scrH}$ taking into account $\X$ and $\supp (\scrH)$. For a random vector $\bV = (V_1, \cdots, V_d)\in\reals^d$ with distribution $\P_V$, the Rosenblatt transformation is
\begin{align*}
	T_{\P_V}(v_1, \cdots, v_d)  =
	\begin{bmatrix}
		\P_V(V_1 \leq v_1)\\
		\P_V(V_2 \leq v_2 | V_1 = v_1)\\
		\vdots \\
		\P_V(V_d \leq v_d| V_1 = v_1, \cdots, V_{d-1} = v_{d-1})
	\end{bmatrix}.
\end{align*}

To discuss the role of the Rosenblatt transformation in constructing a valid normalizing transformation, we consider the following two cases.

\begin{enumerate}
	\item[\textbf{Case 1:}] \textbf{Rectangular $\Q_{\scrH}$}. Suppose that $\bbQ_{\scrH}$ admits a Lebesgue density $f_Q$ supported on $\Q_{\scrH} = \times_{l = 1}^d [\mathsf{a}_l, \mathsf{b}_l]$, $- \infty \leq \mathsf{a}_l < \mathsf{b}_l \leq \infty$. Then, the Rosenblatt transformation $\phi_{\scrH} = T_{\bbQ_\scrH}$ is a normalizing transformation, and we obtain
	\begin{align*}
		\ttc_1 & = d \sup_{\bu \in \Q_{\scrH}} \frac{f_Q(\bu)}{\min \{f_{Q,1}(u_1), f_{Q,2|1}(u_2|u_1), \cdots, f_{Q,d|-d}(u_d|u_1,\cdots,u_{d-1})\}}, \\
		\ttc_2 & = \sup_{\bu \in \Q_{\scrH}} \frac{1}{\min \{f_{Q,1}(u_1), f_{Q,2|1}(u_2|u_1), \cdots, f_{Q,d|-d}(u_d|u_1,\cdots,u_{d-1})\}},
	\end{align*}
	and $\ttc_3 = 2^{d-1} d^{d/2-1} \ttc_1 \ttc_2^{d-1}$, where $f_{Q,j|-j}(\cdot | u_1,\cdots,u_{j-1})$ denotes the conditional density of $Q_j|Q_1=u_1,\cdots,Q_{j-1}=u_{j-1}$ for $\bQ = (Q_1,\cdots,Q_d)\thicksim\bbQ_\scrH$.
	
	This case covers several examples of interest, which give primitive conditions for assumption (ii) in Theorem \ref{thm: X-Process -- Main Theorem}:
	\begin{enumerate}[label=(\alph*)]
		\item Suppose $\Q_\scrH = \times_{l=1}^d [\mathsf{a}_l,\mathsf{b}_l]$ is bounded. Then, for $f_Q$ bounded and bounded away from zero on $\Q_\scrH$,
		\begin{align*}
			\ttc_1 \leq d \frac{\overline{f}_Q^2}{\underline{f}_Q} \overline{\Q}_{\scrH}
			\qquad\text{and}\qquad
			\ttc_2 \leq \frac{\overline{f}_Q}{\underline{f}_Q} \overline{\Q}_{\scrH},
		\end{align*}
		where $\underline{f}_Q = \inf_{\bx \in \Q_\scrH}f_X(\bx)$, $\overline{f}_Q = \sup_{\bx \in \Q_\scrH}f_Q(\bx)$, and $\overline{\Q}_{\scrH} = \max_{1 \leq l \leq d} |\mathsf{b}_l - \mathsf{a}_l|$. 
		If $\X = \times_{l=1}^d [\mathsf{a}_l,\mathsf{b}_l]$ is bounded and $\P_X$ admits a bounded Lebesgue density $f_X$ on $\X$, then we can set $\bbQ_{\scrH} = \P_X$ and $\phi_{\scrH} = T_{\P_X}$. This case corresponds to \citep[Theorem 1.1]{Rio_1994_PTRF}, and the bounds for $\ttc_1$ and $\ttc_3$ coincide with those in \citep[Section 3, \emph{Transformation of the r.v.'s}]{Rio_1994_PTRF}. Alternatively, if $\X$ is unbounded but $\supp(\scrH)$ is bounded, we may still be able to find $\bbQ_\scrH$ supported on a bounded rectangle. We illustrate this case with Example~\ref{example: kde} in Section \ref{sec: Special Cases of Theorem 1}.
		
		\item Suppose $\Q_\scrH = \times_{l = 1}^d [\mathsf{a}_l,\mathsf{b}_l]$ is unbounded. This is often the case when $\X$ and $\supp(\scrH)$ are unbounded (but note that setting $\X \cap \supp(\scrH)$ could be bounded in some cases). To fix ideas, let $\bx_i \thicksim \mathsf{Normal} (\boldsymbol{\mu}, \boldsymbol{\Sigma})$. Then, we can set $\bbQ_{\scrH} = \P_X$ and $\phi_{\scrH} = T_{\P_X}$, and obtain
		\begin{align}
			\ttc_1 & \leq d \sup_{\bx \in \Q_{\scrH}} \max \{f_{X,1}(x_1), f_{X,2|1}(x_2|x_1), \cdots, f_{X,d|-d}(x_d|x_{-d})\}^{d-1} \label{eq: c1-bound X-unbounded}\\
			& \leq d \min_{1 \leq k \leq d} \{\bSigma_{k,k} - \bSigma_{k,1:k-1} \bSigma_{1:k-1,1:k-1}^{-1}\bSigma_{1:k-1,k}\}^{-(d-1)/2}\nonumber
		\end{align}
		bounded, but $\ttc_2$ (and hence $\ttc_3$) unbounded. This result shows that even when the support of $\P_X$ is unbounded, a valid uniform Gaussian strong approximation can be established in certain cases (albeit the Lipschitz property is not used).
	\end{enumerate}
	
	\medskip
	\item[\textbf{Case 2:}] \textbf{Non-Rectangular $\Q_{\scrH}$}. Due to the irregularity of $\X$ and $\supp(\scrH)$, in some settings only a surrogate measure $\bbQ_\scrH$ with non-rectangular $\Q_\scrH$ may exist. Then, we can compose the Rosenblatt transformation with another mapping capturing the shape of $\Q_{\scrH}$ to exhibit a valid normalizing transformation. Suppose that $\bbQ_{\scrH}$ admits a Lebesgue density $f_Q$ supported on $\Q_{\scrH}$, and there exists a diffeomorphism $\chi:\Q_{\scrH}\mapsto[0,1]^d$. Setting $\phi_{\scrH} = T_{\bbQ_\scrH \circ \chi^{-1}} \circ \chi$ gives a valid normalizing transformation, with
	\begin{align*}
		\ttc_1 \leq d\frac{\overline{f}_Q^2}{\underline{f}_Q} \mathtt{S}_\chi
		\qquad\text{and}\qquad
		\ttc_2 \leq \frac{\overline{f}_Q}{\underline{f}_Q} \mathtt{S}_\chi,
	\end{align*}
	where $\mathtt{S}_\chi = \frac{\sup_{\bx \in [0,1]^d} |\operatorname{det}(\nabla\chi^{-1}(\bx))|}{\inf_{\bx \in [0,1]^d} |\operatorname{det}(\nabla\chi^{-1}(\bx))|}\infnorm{\norm{\nabla \chi^{-1}}_2}$. See also Example~\ref{example: kde} in Section \ref{sec: Special Cases of Theorem 1}.    
\end{enumerate}

To recap, Theorem \ref{thm: X-Process -- Main Theorem} requires the existence of a surrogate measure and a normalizing transformation, which restrict the probability law of the data and take advantage of specific features of the function class. In particular, assumption (ii) in Theorem \ref{thm: X-Process -- Main Theorem} does not require $\X$ to be compact if either \eqref{eq: c1-bound X-unbounded} is bounded (as it occurs when $\P_X$ is the Gaussian distribution) or $\supp(\scrH)$ is bounded (as we illustrate in Example~\ref{example: kde} in Section \ref{sec: Special Cases of Theorem 1}). See Section SA-II.2 of \cite{Cattaneo-Yu_2024_AOS--SA} for details.

\subsection{Special Cases and Related Literature}\label{sec: Special Cases of Theorem 1}

We introduce our first statistical example.

\begin{example}[Kernel Density Estimation]\label{example: kde}
	Suppose that $\P_X$ admits a continuous Lebesgue density $f_{X}$ on its support $\X$. The classical kernel density estimator is
	\begin{align*}
		\widehat{f}_X(\bw) = \frac{1}{n} \sum_{i = 1}^n \frac{1}{b^d}K\Big(\frac{\bx_i-\bw}{b}\Big),
	\end{align*}
	where $K: \calK \to \reals$ is a continuous function with $\calK\subseteq\reals^d$ compact, and $\int_{\calK} K(\bw) d\bw = 1$. In statistical applications, the bandwidth $b\to0$ as $n\to\infty$ to enable nonparametric estimation \citep{Wand-Jones_1995_Book}. Consider establishing a strong approximation for the localized empirical process $(\xi_n(\bw):\bw \in\W)$, $\W\subseteq\X$, where
	\begin{equation*}
		\xi_n(\bw)
		= \sqrt{n b^d} \big(\widehat{f}_X(\bw) - \E[\widehat{f}_X(\bw)]\big)
		= X_n(h_{\bw}), \qquad h_{\bw}\in\scrH,
	\end{equation*}
	with $\scrH = \{h_{\bw}(\cdot) = b^{-d/2}K((\cdot - \bw)/b):\bw\in\W\}$. It follows that $\ttM_{\scrH,\reals^d} = O(b^{-d/2}$).
\end{example}

Variants of Example \ref{example: kde} have been discussed extensively in prior literature on strong approximations because the process $\xi_n$ is non-Donsker whenever $b\to0$, and hence standard weak convergence results for empirical processes can not be used. For example, \cite{Gine-Koltchinskii-Sakhanenko_2004_PTRF} and \cite{Gine-Nickl_2010_AoS} established strong approximations for the univariate case ($d=1$) under i.i.d. sampling with $\X$ unbounded, \cite{Cattaneo-Jansson-Ma_2024_JOE} established strong approximations for the univariate case ($d=1$) under i.i.d. sampling with $\X$ compact, \cite{Rio_1994_PTRF} established strong approximations for the multivariate case ($d>1$) under i.i.d. sampling with $\X$ compact, \cite{Sakhanenko_2015_TPA} established strong approximations for the multivariate case ($d>1$) under i.i.d. sampling with $\X$ unbounded, and \cite{Cattaneo-Feng-Underwood_2024_JASA} established strong approximations for the univariate case ($d=1$) under non-i.i.d. dyadic data with $\X$ compact. \cite[Remark 3.1]{chernozhukov2014gaussian} provides further discussion and references. See also \cite{Cattaneo-Chandak-Jansson-Ma_2024_Bernoulli} for an application of \cite{Rio_1994_PTRF} to uniform inference for conditional density estimation.

We can use Example \ref{example: kde} to further illustrate the role of $\bbQ_\scrH$ and $\phi_\scrH$.

\setcounter{example}{0}
\begin{example}[continued]\label{example: kde -- normalizing transformation}
	Recall that $\X$ is the support of $\P_X$, $\W\subseteq\X$ is the index set for the class $\scrH$, and $\calK$ is the compact support of $K$. It follows that $\supp(\scrH) = \W + b \cdot \calK$. We illustrate two sets of primitive conditions implying assumption (ii) in Theorem \ref{thm: X-Process -- Main Theorem}. 
	
	\begin{itemize}
		
		\item Suppose that $\X=\times_{l = 1}^d [\mathsf{a}_l, \mathsf{b}_l]$, $- \infty \leq \mathsf{a}_l < \mathsf{b}_l \leq \infty$, and $\W$ is arbitrary. Then, we can set $\bbQ_{\scrH} = \P_X$ and $\phi_{\scrH} = T_{\P_X}$, and the discussion in parts (a) and (b) of Case 1 in Section~\ref{section: normalizing transformation} applies, which implies assumption (ii) in Theorem \ref{thm: X-Process -- Main Theorem} under the assumptions imposed therein. Furthermore, when $\X$ is bounded, $\ttc_1 = O(1)$ and $\ttc_2 = O(1)$, and hence $\ttc_3 = O(1)$, because $f_X$ is continuous and positive on $\X$. This is part (a) in Case 1 of Section \ref{section: normalizing transformation}, and also the example in \citep[Section 4]{Rio_1994_PTRF}. No information on $\supp(\scrH)$ is used.
		
		\item Suppose that $\X$ is arbitrary, and $\W$ is bounded. Then, it may be possible to find $\bbQ_\scrH$ supported on a bounded set, even if $\X$ is unbounded. For example, suppose that $\X=\reals^d_+$, $\W=\times_{l = 1}^d [\mathsf{a}_l, \mathsf{b}_l]$, $0 \leq \mathsf{a}_l < \mathsf{b}_l < \infty$, and $\calK=[-1,1]^d$. Then, for instance, we can take $\bbQ_\scrH$ with Lebesgue density
		\begin{align*}
			f_Q(\bx) =
			\begin{cases}
				f_X(\bx) & \text{ if } \bx \in \times_{l = 1}^d [\overline{\mathsf{a}}_l, \overline{\mathsf{b}}_l], \\
				(1 - \P_X(\times_{l = 1}^d [\overline{\mathsf{a}}_l, \overline{\mathsf{b}}_l]))/\lebmeas(\Upsilon) & \text{ if } \bx \in \Upsilon, \\
				0 & \text{ otherwise},
			\end{cases}
		\end{align*}
		where $\overline{\mathsf{a}}_l = \max\{\mathsf{a}_l - b,0\}$, $\overline{\mathsf{b}}_l = \mathsf{b}_l + b$, $\Upsilon = \times_{l = 1}^d [\overline{\mathsf{a}}_l, \overline{\mathsf{b}}_l + 1] \setminus \times_{l = 1}^d [\overline{\mathsf{a}}_l, \overline{\mathsf{b}}_l]$, and $\phi_{\scrH} = T_{\bbQ_\scrH \circ \chi^{-1}} \circ \chi$ with $\chi(x_1, \cdots, x_d) = ((\overline{\mathsf{b}}_1 - \overline{\mathsf{a}}_1)^{-1}(x_1 - \overline{\mathsf{a}}_1), \cdots, (\overline{\mathsf{b}}_d - \overline{\mathsf{a}}_d)^{-1}(x_d - \overline{\mathsf{a}}_d))$. It follows that assumption (ii) in Theorem~\ref{thm: X-Process -- Main Theorem} holds. A more general example is discussed in \cite[Section SA-II.6]{Cattaneo-Yu_2024_AOS--SA}.
	\end{itemize}
	Finally, the surrogate measure and normalizing transformation could be used to incorporate truncation arguments. We do not dive into this idea for brevity. 
\end{example}

We now specialize Theorem \ref{thm: X-Process -- Main Theorem} to several cases of practical interest. We employ the definitions and notation conventions given in Section \ref{section:Main Definitions}. To streamline the presentation, we also assume that $\ttc_1 < \infty$ and $\ttc_2 < \infty$ (hence $\ttc_3 < \infty$) in the remaining of Section \ref{section:General Empirical Processes}. See \cite[Section SA-II]{Cattaneo-Yu_2024_AOS--SA} for details.

\subsubsection{VC-type Bounded Functions}\label{section: VC-type Bounded Functions}

Our first corollary considers a VC-type class $\scrH$ of uniformly bounded functions ($\ttM_{\scrH} <\infty$), but without assuming they are Lipschitz ($\ttL_\scrH = \infty$).

\begin{coro}[VC-type Bounded Functions]\label{coro: X-process -- vc bdd}
	Suppose the conditions of Theorem \ref{thm: X-Process -- Main Theorem} hold. In addition, assume that $\scrH$ is a VC-type class with respect to envelope function $\ttM_{\scrH}$ over $\Q_\scrH$ with constants $\ttc_{\scrH} \geq e$ and $\mathtt{d}_{\scrH} \geq 1$. Then, \eqref{eq:SA} holds with
	\begin{equation*}
		\varrho_n
		= \mathsf{m}_{n,d} \sqrt{\log n} \sqrt{\ttc_1 \ttM_{\scrH}  \ttTV_{\scrH}}
		+ \frac{\log n}{\sqrt{n}} \min\{\sqrt{\log n} \sqrt{\ttM_{\scrH}},\sqrt{\ttc_3 \mathtt{K}_{\scrH} + \mathtt{M}_{\scrH}}\} \sqrt{\ttM_{\scrH}}.
	\end{equation*} 
\end{coro}

This corollary recovers the main result in \cite[Theorem 1.1]{Rio_1994_PTRF} when $d\geq2$, where $\mathsf{m}_{n,d}=n^{-1/(2d)}$. It also covers $d=1$, where $\mathsf{m}_{n,1}=n^{-1/2}\sqrt{\log n}$, thereby allowing for a precise comparison with prior KMT strong approximation results in the univariate case \citep{Gine-Koltchinskii-Sakhanenko_2004_PTRF,Gine-Nickl_2010_AoS,Cattaneo-Feng-Underwood_2024_JASA}. Thus, Corollary \ref{coro: X-process -- vc bdd} contributes to the literature by covering all $d \geq 1$ cases simultaneously, allowing for possibly weaker regularity conditions on $\P_X$ through the surrogate measure and normalizing transformation, and making explicit the dependence on $d$, $\X$, and all other features of the underlying data generating process. This additional contribution can be useful for non-asymptotic probability concentration arguments, or for truncation arguments (see \cite{Sakhanenko_2015_TPA} for an example). Nonetheless, for $d\geq2$, the main intellectual content of Corollary \ref{coro: X-process -- vc bdd} is due to \cite{Rio_1994_PTRF}; we present it here for completeness and as a prelude for our upcoming results.

For $d=1$, Corollary \ref{coro: X-process -- vc bdd} delivers the optimal univariate KMT approximation rate when $\mathtt{K}_{\scrH} = O(1)$, which employs a weaker notion of total variation relative to prior literature, but at the expense of requiring additional conditions, as the following remark explains.

\begin{remark}[Univariate Strong Approximation]\label{remark: comparison of tv}
	In Section 2 of \cite{Gine-Koltchinskii-Sakhanenko_2004_PTRF} and the proof of \cite{Gine-Nickl_2010_AoS}, the authors considered univariate ($d=1$) i.i.d. continuously distributed random variables, and established the strong approximation: 
	\begin{align*}
		\P \left(\norm{X_n - Z_n^X}_{\scrH} > \mathtt{pTV}_{\scrH,\reals}\frac{t + C_1 \log n}{\sqrt{n}} \right) \leq C_2 \exp(-C_3 t), \qquad t>0,
	\end{align*}
	where $C_1, C_2, C_3$ are universal constants. \cite[Lemma SA20]{Cattaneo-Feng-Underwood_2024_JASA} slightly generalized the result (e.g., $\P_X$ is not required to be absolutely continuous with respect to the Lebesgue measure), and provided a self-contained proof.
	
	For any interval $\mathcal{I}$ in $\reals$, $\ttTV_{\scrH,\mathcal{I}} \leq \mathtt{pTV}_{\scrH,\mathcal{I}}$ provided that $\ttM_{\scrH,\mathcal{I}} < \infty$ \citep[Theorem 3.27]{Ambrosio-Fusco-Pallara_2000_book}. Therefore, Theorem \ref{thm: X-Process -- Main Theorem} employs a weaker notation of total variation, but imposes complexity requirements on $\scrH$ and the existence of a normalizing transformation. In contrast, \cite{Gine-Koltchinskii-Sakhanenko_2004_PTRF}, \cite{Gine-Nickl_2010_AoS} and \cite{Cattaneo-Feng-Underwood_2024_JASA} do not imposed those extra conditions, but their results only apply when $d = 1$.
\end{remark}

We illustrate the usefulness of Corollary \ref{coro: X-process -- vc bdd} with Example \ref{example: kde}.

\setcounter{example}{0}
\begin{example}[continued]\label{example: kde -- vc bdd}
	Let the conditions of Theorem~\ref{thm: X-Process -- Main Theorem} hold, and $n b^d / \log n \to \infty$. Prior literature further assumed $K$ is Lipschitz to verify the conditions of Corollary~\ref{coro: X-process -- vc bdd} with $\ttTV_{\scrH} = O(b^{d/2-1})$ and $\ttK_{\scrH} = O(b^{-d/2})$. Then, for $X_n=\xi_n$, \eqref{eq:SA} holds with $\varrho_n = (nb^d)^{-1/(2d)} \sqrt{\log n} + (nb^d)^{-1/2} \log n$.
\end{example}

The resulting uniform Gaussian approximation convergence rate in Example \ref{example: kde -- vc bdd} matches prior literature for $d=1$ \citep{Gine-Koltchinskii-Sakhanenko_2004_PTRF,Gine-Nickl_2010_AoS,Cattaneo-Feng-Underwood_2024_JASA} and $d\geq2$ \citep{Rio_1994_PTRF}. This result concerns the uniform Gaussian strong approximation of the entire stochastic process, which can then be specialized to deduce a strong approximation for the scalar suprema of the empirical process $\|\xi_n\|_\scrH$. As noted by \cite[Remark 3.1(ii)]{chernozhukov2014gaussian}, the (almost sure) strong approximation rate in Example \ref{example: kde -- vc bdd} is better than their strong approximation rate (in probability) for $\|\xi_n\|_\scrH$ when $d\in\{1,2,3\}$, but their approach specifically tailored to the scalar suprema delivers better strong approximation rates when $d\geq4$.

Following prior literature, Example \ref{example: kde -- vc bdd} imposed the additional condition that $K$ is Lipschitz to verify that $\scrH = \{b^{-d/2}K((\cdot - \bw)/b):\bw\in \W\}$ forms a VC-type class, and the other conditions in Corollary \ref{coro: X-process -- vc bdd}. The Lipschitz assumption holds for most kernel functions used in practice. One notable exception is the uniform kernel, which is nonetheless covered by Corollary \ref{coro: X-process -- vc bdd}, and prior results in the literature, with a slightly suboptimal strong approximation rate (an extra $\sqrt{\log n}$ term appears when $d\geq2$).

\subsubsection{VC-type Lipschitz Functions}

It is known that the uniform Gaussian strong approximation rate in Corollary~\ref{coro: X-process -- vc bdd} is optimal under the assumptions imposed \citep{beck1985lower}. However, the class of functions $\scrH$ often has additional structure in statistical applications that can be exploited to improve on Corollary~\ref{coro: X-process -- vc bdd}. In Example \ref{example: kde -- vc bdd}, for instance, prior literature further assumed $K$ is Lipschitz to verify the sufficient conditions. Therefore, our next corollary considers a VC-type class $\scrH$ now allowing for the possibility of Lipschitz functions ($\ttL_\scrH < \infty$).

\begin{coro}[VC-type Lipschitz Functions]\label{coro: X-process -- vc lip}
	Suppose the conditions of Theorem \ref{thm: X-Process -- Main Theorem} hold. In addition, assume that $\scrH$ is a VC-type class with envelope function $\ttM_{\scrH}$ over $\Q_\scrH$ with constants $\ttc_{\scrH} \geq e$ and $\mathtt{d}_{\scrH} \geq 1$. Then, \eqref{eq:SA} holds with
	\begin{align*}
		\varrho_n
		&= \min\{\mathsf{m}_{n,d} \sqrt{\ttM_{\scrH}},
		\mathsf{l}_{n,d} \sqrt{\ttc_2 \ttL_{\scrH}}\}  \sqrt{\log n} \sqrt{\ttc_1 \ttTV_{\scrH}} \\
		&\qquad + \frac{\log n}{\sqrt{n}} \min\{\sqrt{\log n} \sqrt{\ttM_{\scrH}},\sqrt{\ttc_3 \mathtt{K}_{\scrH} + \mathtt{M}_{\scrH}}\} \sqrt{\ttM_{\scrH}}.
	\end{align*}
\end{coro}

Putting aside $\ttM_{\scrH}$ and $\ttTV_{\scrH}$, this corollary shows that if $\ttL_\scrH < \infty$, then the rate of strong approximation can be improved. In particular, for $d=2$, $\mathsf{m}_{n,2}=n^{-1/4}$ but $\mathsf{l}_{n,2}=n^{-1/2} \sqrt{\log n}$, implying that $\varrho_n = n^{-1/2} \log n$ whenever $\mathtt{K}_{\scrH} = O(b^{-d/2})$. Therefore, Corollary~\ref{coro: X-process -- vc lip} establishes a uniform Gaussian strong approximation for general empirical processes based on bivariate data that can achieve the optimal univariate KMT approximation rate. (An additional $\sqrt{\log n}$ penalty would appear if $\mathtt{K}_{\scrH}=\infty$.)

For $d\geq3$, Corollary~\ref{coro: X-process -- vc lip} also provides improvements relative to prior literature, but falls short of achieving the optimal univariate KMT approximation rate. Specifically, $\mathsf{m}_{n,d}=n^{-1/(2d)}$ but $\mathsf{l}_{n,d}=n^{-1/d}$ for $d\geq3$, implying that $\varrho_n = n^{-1/d} \sqrt{\log n}$. It remains an open question whether further improvements are possible at this level of generality: the main roadblock underlying the proof strategy is related to the coupling approach based on the Tusn\'ady's inequality for binomial counts, which in turn are generated by the aforementioned mean square approximation of the functions $h\in\scrH$ by local constant functions on carefully chosen partitions of $\Q_\scrH$. Our key observation underlying Corollary~\ref{coro: X-process -- vc lip}, and hence the limitation, is that for Lipschitz functions ($\ttL_\scrH < \infty$) both the projection error arising from the mean square approximation and the KMT coupling error by \cite[Theorem 2.1]{Rio_1994_PTRF} can be improved. However, further improvements for smoother functions appear to necessitate an approximation approach that would not generate dyadic binomial counts, thereby rendering current coupling approaches inapplicable.

We revisit the kernel density estimation example to illustrate the power of Corollary~\ref{coro: X-process -- vc lip}.

\setcounter{example}{0}
\begin{example}[continued]\label{example: kde -- vc lip}
	Under the conditions imposed, $\ttL_\scrH = O(b^{-d/2-1})$, and Corollary~\ref{coro: X-process -- vc lip} implies that, for $X_n=\xi_n$, \eqref{eq:SA} holds with $\varrho_n = (nb^d)^{-1/d} \sqrt{\log n} + (nb^d)^{-1/2} \log n$.
\end{example}

Returning to the discussion of \cite[Remark 3.1(ii)]{chernozhukov2014gaussian}, Example \ref{example: kde -- vc lip} shows that our almost sure strong approximation rate for the entire empirical process is now better than their strong approximation (in probability) rate for the scalar suprema $\|\xi_n\|_\scrH = \sup_{\bw \in \mathcal{W}}|\xi_n(\bw)|$ when $d\leq 6$. On the other hand, their approach delivers a better strong approximation rate in probability for $\|\xi_n\|_\scrH$ when $d\geq7$. Our improvement is obtained without imposing additional assumptions because \cite[Section 4]{Rio_1994_PTRF} already assumed $K$ is Lipschitizian for the verification of the conditions imposed by his strong approximation result (cf. Corollary \ref{coro: X-process -- vc bdd}).

\subsubsection{Polynomial-entropy Functions}

\cite{koltchinskii1994komlos} also considered uniform Gaussian strong approximations for the general empirical process under other notions of entropy for $\scrH$, thereby allowing for more complex classes of functions when compared to \cite{Rio_1994_PTRF}. Furthermore, \cite{koltchinskii1994komlos} employed a Haar approximation condition, which plays a similar role as the total variation and the Lipschitz conditions exploited in our paper. To enable a precise comparison to \cite{koltchinskii1994komlos}, the next corollary considers a class $\scrH$ satisfying a polynomial-entropy condition.

\begin{coro}[Polynomial-entropy Functions]\label{coro: X-process -- poly entropy}
	Suppose the conditions of Theorem \ref{thm: X-Process -- Main Theorem} hold, and that $\scrH$ is a polynomial-entropy class with envelope function $\ttM_{\scrH}$ over $\Q_\scrH$ with constants $\mathtt{a}_{\scrH} >0 $ and $0 < \mathtt{b}_{\scrH}<2$. Then, \eqref{eq:SA} holds as follows:
	\begin{enumerate}[label=\emph{(\roman*)}]
		\item If $\ttL_\scrH\leq\infty$, then
		\begin{align*}
			\varrho_n 
			&= \mathsf{m}_{n,d}\sqrt{\ttc_1 \ttM_{\scrH}  \ttTV_{\scrH}}(\sqrt{\log n}+(\ttc_1 \mathsf{m}_{n,d}^2\ttM_{\scrH}^{-1}\ttTV_{\scrH})^{-\frac{\mathtt{b}_{\scrH}}{4}}) \\
			&\quad + \sqrt{\frac{\ttM_{\scrH}}{n}}\min\{\sqrt{\log n}\sqrt{\ttM_{\scrH}}, \sqrt{\ttc_3 \ttK_{\scrH} + \ttM_{\scrH}}\}(\log n +(\ttc_1 \mathsf{m}_{n,d}^2 \ttM_{\scrH}^{-1} \ttTV_{\scrH})^{-\frac{\mathtt{b}_{\scrH}}{2}}),
		\end{align*}
		
		\item If $\ttL_\scrH<\infty$, then
		\begin{align*}
			\varrho_n 
			&= \mathsf{l}_{n,d}\sqrt{\ttc_1 \ttc_2 \ttL_{\scrH}\ttTV_{\scrH}}(\sqrt{\log n}+(\ttc_1 \ttc_2 \mathsf{l}_{n,d}^2\ttM_{\scrH}^{-2}\ttL_{\scrH}\ttTV_{\scrH})^{-\frac{\mathtt{b}_{\scrH}}{4}}) \\
			&\quad + \sqrt{\frac{\ttM_{\scrH}}{n}}\min\{\sqrt{\log n}\sqrt{\ttM_{\scrH}}, \sqrt{\ttc_3 \ttK_{\scrH} + \ttM_{\scrH}}\}(\log n +(\ttc_1 \ttc_2 \mathsf{l}_{n,d}^2 \ttM_{\scrH}^{-2} \ttL_{\scrH}\ttTV_{\scrH})^{-\frac{\mathtt{b}_{\scrH}}{2}}).
		\end{align*}
	\end{enumerate}
\end{coro}

This corollary reports a simplified version of our result, which corresponds to the best possible bound for the discussion in this section. See \cite[Section SA-II]{Cattaneo-Yu_2024_AOS--SA} for the general case. It is possible to apply Corollary \ref{coro: X-process -- poly entropy} to Example~\ref{example: kde}, although the result is suboptimal relative to the previous results leveraging a VC-type condition.

\setcounter{example}{0}
\begin{example}[continued]\label{example: kde -- pol entropy}
	Under the conditions imposed, for any $0 < \mathtt{b}_{\scrH} < 2$, we can take $\mathtt{a}_{\scrH} = \log(d + 1) + d \mathtt{b}_{\scrH}^{-1}$ so that $\scrH$ is a polynomial-entropy class with constants $(\mathtt{a}_{\scrH}, \mathtt{b}_{\scrH})$. Then, Corollary~\ref{coro: X-process -- poly entropy}\emph{(ii)} implies that, for $X_n=\xi_n$, \eqref{eq:SA} holds with
	$\varrho_n 
	= \mathtt{a}_{\scrH}^2 (n b^d)^{-\frac{1}{d}(1 - \frac{\mathtt{b}_{\scrH}}{2})} b^{-d \mathtt{b}_{\scrH}} 
	+ \mathtt{a}_{\scrH}^2 (n b^d)^{-\frac{1}{2}+\frac{\mathtt{b}_{\scrH}}{d}} b^{-\frac{d \mathtt{b}_{\scrH}}{2}}
	$.
\end{example}

Our running example shows that a uniform Gaussian strong approximation based on polynomial-entropy conditions can lead to suboptimal KMT approximation rates. However, for other (larger) function classes, those results may be useful. The following remark discusses an example studied in \cite{koltchinskii1994komlos}, and illustrates our contributions in that context.

\begin{remark}[Polynomial-entropy Condition]\label{remark: C^q -- pol entropy}
	Suppose $\P_X$ is $\mathsf{Uniform}(\X)$ with $\X=[0,1]^d$, and $\scrH$ a subclass of $C^q(\X)$ with $C^q$-norm uniformly bounded by $1$ and $2 \leq d < q$. \cite[page 111]{koltchinskii1994komlos} discusses this example after his Theorem 11.3, and reports the uniform Gaussian strong approximation rate $n^{-\frac{q - d}{2qd}} \polylog(n)$. See \cite{koltchinskii1994komlos}, or \cite[Section SA-I]{Cattaneo-Yu_2024_AOS--SA}, for the additional notation and definitions used in this example.
	
	Corollary~\ref{coro: X-process -- poly entropy} is applicable to this case, upon setting $(\bbQ_\scrH,\phi_\scrH)=(\P_X, \operatorname{Id})$ with $\operatorname{Id}$ denoting the identity map from $[0,1]^d$ to $[0,1]^d$. It follows that $\mathtt{M}_{\scrH} = 1$, $\mathtt{TV}_{\scrH} = 1$, $\mathtt{L}_{\scrH} = 1$. \cite[Theorem 2.7.1]{wellner2013weak} shows that $\scrH$ is a polynomial-entropy class with constants $\mathtt{a}_{\scrH} = \mathtt{C}_{q,d}$ and $\mathtt{b}_{\scrH}= d/q$, where $\mathtt{C}_{q,d}$ is a constant depending on $q$ and $d$ only. Then, Corollary~\ref{coro: X-process -- poly entropy}\emph{(ii)} implies that, for $X_n=\xi_n$, \eqref{eq:SA} holds with
	\begin{align*}  
		\varrho_n =  
		\begin{cases}
			n^{-\frac{1}{2}+\frac{1}{q}} \polylog(n)  & \text{ if } d = 2\\
			n^{-\frac{2 q - d}{2 d q}} \polylog(n)    & \text{ if } d > 2
		\end{cases},
	\end{align*}
	which gives a faster convergence rate than the one obtained by \cite{koltchinskii1994komlos}.
	
	The improvement is explained by two differences between \cite{koltchinskii1994komlos} and our approach. First, we explicitly incorporate the Lipschitz condition, and hence we can take $\beta = \frac{2}{d}$ instead of $\beta = \frac{1}{d}$ in Equation (3.1) of \cite{koltchinskii1994komlos}. Second, using the uniform entropy condition approach, we get $\log N(\mathscr{H},\norm{\cdot}_{\P_X,2}, \varepsilon) = O(\varepsilon^{-d/q}$), while \cite{koltchinskii1994komlos} started with the bracketing number condition $\log N_{[\;]}(\mathscr{H}, \norm{\cdot}_{\P_X,1}, \varepsilon) = O(\varepsilon^{-d/q})$ and, with the help of his Lemma 8.4, applied Theorem 3.1 with $\alpha = \frac{d}{d + q}$ in his Equation (3.2). The proof of his Theorem 3.1 leverages the fact that his Equation (3.2) implies that $\log N(\mathscr{H}, \norm{\cdot}_{\P_X,2}, \varepsilon) = O(\varepsilon^{-2d/q})$, and his approximation rate is looser by a power of two when compared to the uniform entropy condition underlying our Corollary~\ref{coro: X-process -- poly entropy}. Setting $\mathtt{L}_{\scrH} = \infty$, $\mathtt{b}_{\scrH} = 2d/q$, and keeping the other constants, Corollary~\ref{coro: X-process -- poly entropy}\emph{(i)} would give $\varrho_n = n^{-\frac{q - d}{2 q d}} \polylog(n)$, which is the same rate as in \cite{koltchinskii1994komlos}. Finally, Theorem 3.2 in \cite{koltchinskii1994komlos} allows for $\log N(\mathscr{H}, \norm{\cdot}_{\P_X,2}, \varepsilon) = O(\varepsilon^{-2\rho})$ where $\rho$ is not implied by his Equation (3.2), and his result would give the strong approximation rate $n^{-\frac{2 q - d}{4 qd}} \polylog(n)$.
\end{remark}

\section{Residual-Based Empirical Process}\label{section: Residual Empirical Process}

Consider the simple local empirical process discussed in \cite[Section 3.1]{chernozhukov2014gaussian}:
\begin{align}\label{eq: local empirical process}
	S_n(\bw) = \frac{1}{nb^d} \sum_{i=1}^n K\Big(\frac{\bx_i-\bw}{b}\Big) y_i, \qquad \bw\in\W,
\end{align}
where $\bx_i\thicksim\P_X$, $y_i\thicksim \P_Y$, and $b\to0$ as $n\to\infty$. Using our notation, $\big(\sqrt{n b^d}(S_n(\bw) - \E[S_n(\bw)|\bx_1, \cdots, \bx_n]): \bw \in \W\big) = (R_n(g,r): g \in \G, r \in \R)$ with $\G = \{b^{-d/2}K(\frac{\cdot - \bw}{b}):\bw\in \W\}$ and $\R=\{\operatorname{Id}\}$, where $\operatorname{Id}$ denotes the identity map from $\reals$ to $\reals$. This setting corresponds to kernel regression estimation with $K$ interpreted as the equivalent kernel; see Section \ref{section: Local Polynomial Regression} for details. As noted in \cite[Remark 3.1(iii)]{chernozhukov2014gaussian}, a direct application of \cite{Rio_1994_PTRF}, or of our Theorem \ref{thm: X-Process -- Main Theorem}, views $\bz_i=(\bx_i,y_i)\thicksim\P_Z$ as the underlying $(d+1)$-dimensional random vectors entering the general empirical process $X_n$ defined in \eqref{eq:X-process}. Specifically, under some regularity conditions on $K$ and non-trivial restrictions on the joint distribution $\P_Z$, \cite{Rio_1994_PTRF}'s strong approximation result verifies \eqref{eq:SA} with rate \eqref{eq: local empirical process -- Rate Rio}, which is also verified via Corollary \ref{coro: X-process -- vc bdd}. Furthermore, imposing a Lipschitz property on $\scrH=\G\times\R$, Corollary \ref{coro: X-process -- vc lip} would give the improved strong approximation result \eqref{eq: local empirical process -- Rate Thm 1}, under regularity conditions.

The strong approximation results for $S_n$ illustrate two fundamental limitations because all the elements in $\bz_i=(\bx_i,y_i)$ are treated symmetrically. First, the effective sample size emerging in the strong approximation rate is $nb^{d+1}$, which is suboptimal because only the $d$-dimensional covariate $\bx_i$ are being smoothed out. Since the pointwise variance of the process is of order $n^{-1}b^{-d}$, the correct effective sample size should be $nb^{d}$, up to $\polylog(n)$ terms. Therefore, applying \cite{Rio_1994_PTRF}, or our improved Theorem \ref{thm: X-Process -- Main Theorem}, leads to a suboptimal uniform Gaussian strong approximation for $S_n$. Second, applying \cite{Rio_1994_PTRF}, or our improved Theorem \ref{thm: X-Process -- Main Theorem}, requires $\P_Z$ to be continuously distributed and supported on $[0,1]^{d+1}$, possibly after applying a normalizing transformation. This requirement imposes non-trivial restrictions on $\P_Z$ and, in particular, on $\P_Y$, limiting the applicability of the strong approximation results. See \cite[Remark 3.1(iii)]{chernozhukov2014gaussian} for more discussion.

Motivated by the aforementioned limitations, the following theorem explicitly studies the residual-based empirical process defined in \eqref{eq:R-process}, leveraging its intrinsic multiplicative separable structure. We present our result under a VC-type condition on $\G\times\R$ to streamline the discussion, but a result at the same level of generality as Theorem \ref{thm: X-Process -- Main Theorem} is given in \cite[Section SA-IV]{Cattaneo-Yu_2024_AOS--SA}. Recall Section \ref{section:Main Definitions} and the notation conventions introduced therein.

\begin{thm}\label{thm: R-Process -- Main Theorem}
	Suppose $(\bz_i=(\bx_i, y_i): 1 \leq i \leq n)$ are i.i.d. random vectors taking values in $(\reals^{d+1}, \mathcal{B}(\reals^{d+1}))$ with common law $\P_Z$, where $\bx_i$ has distribution $\P_X$ supported on $\X\subseteq\reals^d$, $y_i$ has distribution $\P_Y$ supported on $\Y\subseteq\reals$, and the following conditions hold.
	\begin{enumerate}[label=(\roman*)]
		\item $\G$ is a real-valued pointwise measurable class of functions on $(\reals^d, \mathcal{B}(\reals^d), \P_X)$.
		\item There exists a surrogate measure $\bbQ_\G$ for $\P_X$ with respect to $\G$ such that $\bbQ_\G = \lebmeas \circ \phi_\G$, where the \textit{normalizing transformation} $\phi_{\G}: \Q_\G \mapsto [0,1]^d$ is a diffeomorphism.
		\item $\G$ is a VC-type class with function $\ttM_{\G}$ over $\Q_\G$ with $\ttc_{\G} \geq e$ and $\mathtt{d}_{\G} \geq 1$.        
		\item $\R$ is a real-valued pointwise measurable class of functions on $(\reals, \mathcal{B}(\reals),\P_Y)$.
		\item $\R$ is a VC-type class with envelope $M_{\R,\Y}$ over $\Y$ with $\ttc_{\R,\Y}\geq e$ and $\mathtt{d}_{\R,\Y}\geq 1$, where $M_{\R,\Y}(y) + \mathtt{pTV}_{\R,(-|y|,|y|)} \leq \ttv (1 + |y|^{\alpha})$ for all $y \in \Y$, for some $\ttv>0$, and for some $\alpha\geq0$. Furthermore, if $\alpha>0$, then $\sup_{\bx \in \X}\E[\exp(|y_i|)|\bx_i = \bx] \leq 2$.
		\item There exists a constant $\ttk$ such that $|\log_2 \mathtt{E}_{\G}| + |\log_2 \ttTV| + |\log_2 \ttM_{\G}| \leq \ttk \log_2 n$, where $\ttTV = \max \{\ttTV_{\G}, \ttTV_{\G \times \mathscr{V}_{\R},\Q_\G}\}$ with $\mathscr{V}_{\R} = \{\theta(\cdot,r) : r \in \R\}$, and $\theta(\bx,r) = \E[r(y_i)|\bx_i = \bx]$.
	\end{enumerate}
	Then, on a possibly enlarged probability space, there exists a sequence of mean-zero Gaussian processes $(Z_n^R(g,r): (g,r)\in \G\times \R)$ with almost sure continuous trajectories on $(\G \times \R, \d_{\P_Z})$ such that:
	\begin{itemize}
		\item $\E[R_n(g_1, r_1) R_n(g_2, r_2)] = \E[Z^R_n(g_1, r_1) Z^R_n(g_2, r_2)]$ for all $(g_1, r_1), (g_2, r_2) \in \G \times \R$, and
		\item $\P\big[\norm{R_n - Z_n^R}_{\G \times \R} > C_1C_{\ttv,\alpha} \mathsf{T}_n(t)\big] \leq C_2 e^{-t}$ for all $t > 0$,
	\end{itemize}
	where $C_1$ and $C_2$ are universal constants, $C_{\ttv,\alpha} = \ttv \max\{1 + (2 \alpha)^{\frac{\alpha}{2}}, 1 + (4 \alpha)^{\alpha}\}$, and
	\begin{gather*}
		\mathsf{T}_n(t)
		=  \mathsf{A}_{n} (t + \ttk \log_2 n + \mathtt{d}\log (\ttc n))^{\alpha + \frac{3}{2}} \sqrt{d}
		+ \frac{\ttM_{\G}}{\sqrt{n}} (t + \ttk \log_2 n+\mathtt{d}\log(\ttc n))^{\alpha + 1}, \\
		\mathsf{A}_{n}
		= \min\bigg\{\bigg( \frac{\ttc_1^d \ttM_{\G}^{d+1} \ttTV^d \mathtt{E}_{\G} }{n}\bigg)^{\frac{1}{2d+2}}, \bigg( \frac{\ttc_1^{\frac{d}{2}} \ttc_2^{\frac{d}{2}} \ttM_{\G} \mathtt{E}_{\G} \ttTV^{\frac{d}{2}} \ttL^{\frac{d}{2}}}{n} \bigg)^{\frac{1}{d+2}} \bigg\},
	\end{gather*}
	\begin{align*}
		\ttc_1 = d \sup_{\bx \in \Q_{\G}} \prod_{j = 1}^{d-1} \sigma_j(\nabla \phi_{\G}(\bx)), \qquad 
		\ttc_2 = \sup_{\bx \in \Q_{\G}} \frac{1}{\sigma_{d}(\nabla \phi_{\G}(\bx))}, 
	\end{align*}
	with $\ttc = \ttc_{\G} \ttc_{\R,\Y}$, $\mathtt{d} = \mathtt{d}_{\G} + \mathtt{d}_{\R,\Y}$, and $\ttL = \max\{\ttL_{\G}$, $\ttL_{\G \times \mathscr{V}_{\R},\Q_\G}\}$.
\end{thm}

This theorem establishes a uniform Gaussian strong approximation under regularity conditions specifically tailored to leverage the multiplicative separable structure of $R_n$ defined in \eqref{eq:R-process}. Conditions (i)--(iii) in Theorem \ref{thm: R-Process -- Main Theorem} are analogous to the conditions imposed in Corollaries \ref{coro: X-process -- vc bdd} and \ref{coro: X-process -- vc lip} for the general empirical process. Conditions (iv)--(v) in Theorem \ref{thm: R-Process -- Main Theorem} are new, mild restrictions on the portion of the stochastic process corresponding to the outcome $y_i$. Condition (v) either assumes $\R$ is uniformly bounded, or restricts the tail decay of the function class $\R$, without imposing restrictive assumptions on the distribution $\P_Y$. Finally, condition (vi) is imposed only to simplify the exposition; see \cite{Cattaneo-Yu_2024_AOS--SA} for the general result. We require a $\mathtt{pTV}$ condition on $\R$ in (v), but $\mathtt{TV}$ conditions on $\G$ and $\G \times \mathscr{V}_{\R}$ in (vi), because $\P_X$ admits a Lebesgue density, but $\P_Y$ may not.

The proof strategy of Theorem \ref{thm: R-Process -- Main Theorem} is similar to the proof for the general empirical process (Theorem \ref{thm: X-Process -- Main Theorem}), and is given in \cite[Section SA-IV]{Cattaneo-Yu_2024_AOS--SA}. First, we discretize to a $\delta$-net to obtain
\begin{align*}
	\norm{R_n - Z_n^R}_{\G \times \R}
	&\leq \norm{R_n - R_n\circ\pi_{(\G \times \R)_\delta}}_{\G \times \R}
	+ \norm{R_n - Z_n^R}_{(\G \times \R)_\delta} \\
	&\qquad + \norm{Z_n^R\circ\pi_{(\G \times \R)_\delta}-Z_n^R}_{\G \times \R},
\end{align*}
where the terms capturing fluctuation off-the-net, $\norm{R_n - R_n\circ\pi_{(\G \times \R)_\delta}}_{\G \times \R}$ and $\norm{Z_n^R\circ\pi_{(\G \times \R)_\delta}-Z_n^R}_{\G \times \R}$, are handled via standard empirical process methods. Second, the remaining term $\norm{R_n - Z_n^R}_{(\G \times \R)_\delta}$, which captures the finite-class Gaussian approximation error, is once again decomposed via a suitable mean square projection onto the class of piecewise constant Haar functions on a carefully chosen collection of cells partitioning the support of $\P_Z$. This is our point of departure from prior literature.

We design the partitioning cells based on two key observations: (i) regularity conditions are often imposed on the conditional distribution of $y_i|\bx_i$, as opposed to on their joint distribution; and (ii) $\G$ and $\R$ often require different regularity conditions. For example, in the classical regression case discussed previously, $\R$ is just the singleton identity function but $\P_Y$ may have unbounded support or atoms, while $\G$ is a VC-type class of $n$-varying functions with a possibly more regular $\P_X$ having compact support. Furthermore, the dimension of $y_i$ is a nuisance for the strong approximation, making results like Theorem \ref{thm: X-Process -- Main Theorem} suboptimal in general. These observations suggest choosing dyadic cells by an asymmetric iterative splitting construction, where first the support of each dimension of $\bx_i$ is partitioned, and only after the support of $y_i$ is partitioned based on the conditional distribution of $y_i|\bx_i$. See \cite{Cattaneo-Yu_2024_AOS--SA} for details on our proposed asymmetric dyadic cells expansion.

Given our dyadic expansion exploiting the structure of $R_n$, we decompose the term $\norm{R_n - Z_n^R}_{(\G \times \R)_\delta}$ similarly to \eqref{eq: X-Process -- Coupling}, leading to a projected piecewise constant process and the corresponding two projection errors. However, instead of employing the $L_2$-projection $\proj$ as in \eqref{eq: X-Process -- Coupling}, we now use another mapping $\projreg$ from $L_2(\P_Z)$ to piecewise constant functions that explicitly factorizes the product $g(\bx_i)r(y_i)$. In fact, as we discuss in \cite{Cattaneo-Yu_2024_AOS--SA}, each base level cell $\C$ produced by our asymmetric dyadic splitting scheme can be written as a product of the form $\X_l \times \Y_m$, where $\X_l$ denotes the $l$-th cell for $\bx_i$ and $\Y_m$ denotes the $m$-th cell for $y_i$. Thus, $\projreg$ is carefully chosen so that once we know $\bx \in \X_l$ for some $l$, $\projreg[g,r](\bx,y) = \sum_{m=0}^{2^N - 1} \Indicator(y \in \Y_m)\E[r(y_i)|y_i \in \Y_m, \bx_i \in \X_l]\E[g(\bx_i)|\bx_i \in \X_l]$, which only depends on $y$, and has envelope and total variation no greater than those for $r$.

Finally, our generalized Tusn\'ady's lemma for more general binomial counts \citep{Cattaneo-Yu_2024_AOS--SA} allows for the Gaussian coupling of any piecewise-constant functions over our asymmetrically constructed dyadic cells. A generalization of \cite[Theorem 2.1]{Rio_1994_PTRF} enables upper bounding the Gaussian approximation error for processes indexed by piecewise constant functions by summing up a quadratic variation from all layers in the cell expansion. By the above choice of cells and projections, the contribution from the last layers corresponding to splitting $y_i$ amounts to a sum of one-dimensional KMT coupling error from all possible $\X_l$ cells. In fact, the one-dimensional KMT coupling is optimal and, as a consequence, requiring a vanishing contribution of $y_i$ layers to the approximation error does not add extra requirements besides conditions on envelope functions and an $L_1$ bound for $\G$. This explains why we can obtain strong approximation rates reflecting the correct effective sample size underlying the empirical process for the kernel regression and other local empirical process examples.

The following corollary summarizes the main result from Theorem \ref{thm: R-Process -- Main Theorem}.

\begin{coro}[VC-Type Lipschitz Functions]\label{coro: R-process -- vc lip}
	Suppose the conditions of Theorem \ref{thm: R-Process -- Main Theorem} hold with constants $\ttc$ and $\ttd$. Then, $\norm{R_n - Z_n^R}_{\G \times \R} = O(\varrho_n)$ a.s. with
	\begin{equation*}
		\varrho_n
		= \min\Big\{\frac{(\ttc_1^d \ttM_{\G}^{d+1} \ttTV^d \mathtt{E}_{\G})^{\frac{1}{2d+2}} }{n^{1/(2d+2)}},
		\frac{(\ttc_1^{\frac{d}{2}}\ttc_2^{\frac{d}{2}}\ttM_{\G} \ttTV^{\frac{d}{2}} \mathtt{E}_{\G} \ttL^{\frac{d}{2}})^{\frac{1}{d+2}}}{n^{1/(d+2)}} \Big\} (\log n)^{\alpha+3/2}
		+ \frac{(\log n)^{\alpha+1}}{\sqrt{n}} \ttM_{\G}.
	\end{equation*}
\end{coro}

This corollary shows that our best attainable uniform Gaussian strong approximation rate for $R_n$ is $n^{-1/(d+2)} \polylog(n)$, putting aside $\ttc_1$, $\ttc_2$, $\ttM_{\G}$, $\ttTV$, $\mathtt{E}_{\G}$, and $\ttL$. It is not possible to give a strict ranking between Corollary \ref{coro: X-process -- vc lip} and Corollary \ref{coro: R-process -- vc lip}. On the one hand, Corollary \ref{coro: X-process -- vc lip} treats all components in $\bz_i$ symmetrically, and thus imposes stronger regularity conditions on $\P_Z$, but leads to the better approximation rate $n^{-\min\{1/(d+1),1/2\}} \polylog(n)$, putting aside the various constants and underlying assumptions. On the other hand, Corollary \ref{coro: R-process -- vc lip} can deliver a tighter strong approximation under weaker regularity conditions whenever $\scrH = \G \times \R$  and $\G $ varies with $n$, as in the case of the local empirical processes arising from nonparametric regression. The next section offers an application illustrating this point.

See \citep[Section SA-IV]{Cattaneo-Yu_2024_AOS--SA} for proofs and other omitted details. In addition, Section SA-III in \citep{Cattaneo-Yu_2024_AOS--SA} present uniform Gaussian strong approximation results for a general multiplicative-separable empirical process, which may be of interest but is not discussed in the paper to conserve space.

\subsection{Example: Local Polynomial Regression}\label{section: Local Polynomial Regression}

Suppose that $(\bx_1, y_1),\dots,(\bx_n, y_n)$ are i.i.d random vectors taking values in $(\reals^{d+1}, \mathcal{B}(\reals^{d+1}))$, with $\bx_i\thicksim\P_X$ admitting a continuous Lebesgue density on its support $\X = [0,1]^d$. Consider the class of estimands
\begin{align}\label{eq: regression estimand}
	\theta(\bw;r) = \E[r(y_i)|\bx_i=\bw], \qquad \bw \in \W\subseteq\X, \quad r\in\R,
\end{align}
where we focus on two leading cases to streamline the discussion: $\R_1 =\{\operatorname{Id}\}$ corresponds to the conditional expectation $\mu(\bw) = \E[y_i|\bx_i=\bw]$, and $\R_2 =\{\Indicator( \cdot \leq y): y\in\reals\}$ corresponds to the conditional distribution function $F(y|\bw) = \E[\Indicator(y_i \leq y)|\bx_i=\bw]$. In the first case, $\R$ is a singleton but the identity function calls for the possibility of $\P_Y$ not being dominated by the Lebesgue measure or perhaps being continuously distributed with unbounded support. In the second case, $\R$ is a VC-type class of indicator functions, and hence $r(y_i)$ is uniformly bounded, but establishing uniformity over $\R$ is of statistical interest (e.g., to construct specification hypothesis tests based on conditional distribution functions).

Suppose the kernel function $K: \reals^d \to \reals$ is non-negative, Lipschitz, and has compact support $\mathcal{K}$. Using standard multi-index notation, $\bp(\bu)$ denotes the $\frac{(d+\pOrder)!}{d!\pOrder!}$-dimensional vector collecting the ordered elements $\bu^{\bnu}/\bnu!$ for $0\leq|\bnu|\leq \pOrder$, where $\bu^{\bnu}=u_1^{\nu_1}\cdots u_d^{\nu_d}$, $\bnu!=\nu_1! \cdots \nu_d!$ and $|\bnu|=\nu_1+\cdots+\nu_d$, for $\bu=(u_1,\cdots,u_d)^\top$ and $\bnu=(\nu_1,\cdots,\nu_d)^\top$. A local polynomial regression estimator of $\theta(\bw;r)$ is
\begin{align*}
	\widehat{\theta}(\bw;r) = \be_1^{\top}\widehat{\bbeta}(\bw,r),\qquad 
	\widehat{\bbeta}(\bw,r) = \argmin_{\bbeta} \sum_{i = 1}^{n} \big(r(y_i) - \bp(\bx_i - \bw)^{\top}\bbeta\big)^2 K\Big(\frac{\bx_i - \bw}{b}\Big),
\end{align*}
with $\bw \in \W\subseteq\X$, $r\in\R_1$ or $r \in \R_2$, and $\be_1$ denoting the first standard basis vector. See \cite{Fan-Gijbels_1996_Book} for a textbook review. The estimation error can be decomposed into three terms:
\begin{equation*}
	\widehat{\theta}(\bw,r) - \theta(\bw,r)
	= \underbrace{ \be_1^{\top} \bH_{\bw}^{-1}\bS_{\bw,r} }_\text{linearization}
	+ \underbrace{ \be_1^{\top} (\widehat{\bH}_{\bw}^{-1} - \bH_{\bw}^{-1}) \bS_{\bw,r} }_\text{non-linearity error}
	+ \underbrace{ \E[\widehat{\theta}(\bw,r)|\bx_1,\cdots,\bx_n] - \theta(\bw,r) }_\text{smoothing bias},
\end{equation*}
with $\widehat{\bH}_{\bw} = \frac{1}{n}\sum_{i=1}^n \bp(\frac{\bx_i - \bw}{b}) \bp(\frac{\bx_i - \bw}{b})^{\top} \frac{1}{b^d} K (\frac{\bx_i - \bw}{b})$, $\bH_{\bw} = \E [\bp(\frac{\bx_i - \bw}{b}) \bp(\frac{\bx_i - \bw}{b})^{\top} \frac{1}{b^d} K (\frac{\bx_i - \bw}{b})]$, and $\mathbf{S}_{\bw,r} = \frac{1}{n}\sum_{i=1}^n \bp(\frac{\bx_i - \bw}{b}) \frac{1}{b^d} K(\frac{\bx_i - \bw}{b})(r(y_i)-\E[r(y_i)|\bx_i])$.

It follows that the linear term is
\begin{align*}
	\sqrt{n b^d}\be_1^{\top} \bH_{\bw}^{-1}\mathbf{S}_{\bw,r}
	= \frac{1}{\sqrt{n b^d}} \sum_{i=1}^n \mathfrak{K}_\bw\Big(\frac{\bx_i-\bw}{b}\Big) (r(y_i)-\E[r(y_i)|\bx_i])
	= R_n(g,r),
\end{align*}
for $(g,r) \in \G\times\R_l$, $l=1,2$, and where $\G = \{b^{-d/2}\mathfrak{K}_\bw(\frac{\cdot - \bw}{b}):\bw\in \W\}$ with $\mathfrak{K}_\bw(\bu)=\be_1^{\top} \bH_{\bw}^{-1} \bp(\bu)K(\bu)$ the equivalent boundary-adaptive kernel function. Furthermore, under the regularity conditions given in \cite[Section SA-IV.6]{Cattaneo-Yu_2024_AOS--SA}, which relate to uniform smoothness and moment restrictions for the conditional distribution of $y_i|\bx_i$,
\begin{align*}
	&\sup_{\bw\in\W,r\in\R_1} \big|\be_1^{\top} (\widehat{\bH}_{\bw}^{-1} - \bH_{\bw}^{-1}) \bS_{\bw,r}\big|
	= O((n b^d)^{-1} \log n + (n b^d)^{-3/2}(\log n)^{5/2}) \quad \text{a.s.},\\
	&\sup_{\bw\in\W,r\in\R_2} \big|\be_1^{\top} (\widehat{\bH}_{\bw}^{-1} - \bH_{\bw}^{-1}) \bS_{\bw,r}\big|
	= O((n b^d)^{-1}\log{n}) \quad \text{a.s.},\\
	&\sup_{\bw\in\W,r\in\R_l} \big| \E[\widehat{\theta}(\bw,r)|\bx_1,\cdots,\bx_n] - \theta(\bw,r)\big|
	= O(b^{1+\pOrder}) \qquad \text{a.s.}, \quad l=1,2,
\end{align*}
provided that $\log(n)/(n b^d)\to0$. Therefore, the goal reduces to establishing a Gaussian strong approximation for the residual-based empirical process $(R_n(g,r): (g,r) \in \G\times \R_l)$, $l=1,2$. We discuss different attempts to establish such approximation result, culminating with the application of our Theorem \ref{thm: R-Process -- Main Theorem}.

As discussed in \cite[][Remark 3.1]{chernozhukov2014gaussian}, a first attempt is to deploy Theorem 1.1 in \cite{Rio_1994_PTRF} (or, equivalently, Corollary \ref{coro: X-process -- vc bdd}). Viewing the empirical process as based on the random sample $\bz_i=(\bx_i,y_i)\thicksim\P_Z$, $i=1,2,\cdots,n$, Theorem 1.1 in \cite{Rio_1994_PTRF} requires $\P_Z$ to be continuously distributed with positive Lebesgue density on its support $[0,1]^{d+1}$. For this reason, \cite[][Remark 3.1]{chernozhukov2014gaussian} assumes that $(\bx_i, y_i) = (\bx_i, \varphi(\bx_i, u_i))$ where the joint law $\P_B$ of $\mathbf{b}_i = (\bx_i, u_i)$ admits a continuous Lebesgue density supported on $\mathcal{B} = [0,1]^{d+1}$. If $\ttM_{\{\varphi\},\mathcal{B}} < \infty$, $\ttK_{\{\varphi\},\mathcal{B}} < \infty$, $\sup_{g \in \G}\mathtt{TV}_{\{\varphi\},\operatorname{supp}(g) \times [0,1]} < \infty$, and other regularity conditions hold, then it can be shown \citep[Section SA-IV.6]{Cattaneo-Yu_2024_AOS--SA} that applying \cite{Rio_1994_PTRF} to $(X_n(h): h \in \scrH_l)$ based on $(\mathbf{b}_i: 1 \leq i \leq n)$ with $\scrH_l = \{g \cdot (r \circ \varphi) - g \cdot \theta(\cdot,r): g \in \G, r \in \R_l\}$, $l=1,2$, gives a Gaussian strong approximation with rate \eqref{eq: local empirical process -- Rate Rio}. Without the local total variation condition $\ttK_{\{\varphi\},\mathcal{B}} < \infty$, an additional $\sqrt{\log n}$ multiplicative factor appears in the final rate.

The previous result does not exploit Lipschitz continuity, so a natural second attempt is to employ Corollary \ref{coro: X-process -- vc lip} to improve it. Retaining the same assumptions, but now also assuming that $\varphi$ is Lipschitz, our Theorem~\ref{thm: X-Process -- Main Theorem} gives a Gaussian strong approximation for $(X_n(h): h \in \scrH_1)$ with rate \eqref{eq: local empirical process -- Rate Thm 1}. Theorem \ref{thm: X-Process -- Main Theorem} does not give an improvement for $\R_2$ because the Lipschitz condition is not satisfied. See \cite[Section SA-IV.6]{Cattaneo-Yu_2024_AOS--SA}.

The two attempts so far impose restrictive assumptions on the joint distribution of the data, and deliver approximation rates based on the incorrect effective sample size (and thus require $n b^{d+1} \to \infty$). Our Theorem \ref{thm: R-Process -- Main Theorem} addresses both problems: since $\supp(\scrH) = \W + b \mathcal{K}$, and under standard regularity conditions, we can set $\bbQ_\scrH$ and $\phi_\scrH$ according to the discussion in Example~\ref{example: kde -- normalizing transformation}, and thus we verify in \cite[Section SA-IV.6]{Cattaneo-Yu_2024_AOS--SA} that $\ttc_1=O(1)$, $\ttc_2=O(1)$, $\ttM_{\G} = O(b^{-d/2})$, $\ttE_{\G} = O(b^{d/2})$, $\ttK_{\G} = O(b^{-d/2})$, $\ttTV = O(b^{d/2-1})$, and $\ttL = O(b^{-d/2-1})$. This gives $\norm{R_n - Z_n^R}_{\G \times \R_2} = O(\varrho_n)$ a.s. with
\begin{align*}
	\varrho_n = (n b^d)^{-1/(d+2)} \sqrt{\log n} + (n b^d)^{-1/2} \log n.
\end{align*}
If, in addition, we assume $\sup_{\bw \in \W}\E[\exp(|y_i|)|\bx_i = \bw] < \infty$, then $\norm{R_n - Z_n^R}_{\G \times \R_1} = O(\varrho_n)$ a.s. with
\begin{align*}
	\varrho_n = (n b^d)^{-1/(d+2)} \sqrt{\log n} + (n b^d)^{-1/2} (\log n)^2.
\end{align*}

As a consequence, our results verify that the following strong approximations hold:
\begin{itemize}
	\item Let $\widehat{\mu}(\bw)= \widehat{\theta}(\bw;r)$ for $r\in\R_1$. Recall that $\R_1$ consists of the singleton of identity function $\operatorname{Id}$. If $b^{\pOrder + 1}(n b^d)^{\frac{d+4}{2d+4}}(\log n)^{-1/2} + (n b^d)^{-\frac{d+1}{d+2}} (\log n)^2 = O(1)$, then
	\begin{align*}
		\sup_{\bw\in\W} \big|\sqrt{n b^d}\big(\widehat{\mu}(\bw) - \mu(\bw)\big) - Z_n^R(\bw)\big|
		= O(\mathsf{r}_n) \quad \text{a.s.},
		\qquad 
		\mathsf{r}_n = \Big(\frac{(\log n)^{1+d/2}}{n b^d}\Big)^{\frac{1}{d+2}},
	\end{align*}
	where $\cov(Z_n^R(\bw_1),Z_n^R(\bw_2)) = n b^d \cov(\be_1^{\top}\bH_{\bw_1}^{-1}\mathbf{S}_{\bw_1,\operatorname{Id}}, \be_1^{\top}\bH_{\bw_2}^{-1}\mathbf{S}_{\bw_2,\operatorname{Id}})$ for all $\bw_1,\bw_2\in\W$.
	
	\item Let $\widehat{F}(y|\bw)= \widehat{\theta}(\bw;r_{y})$ for $r_{y} = \Indicator(\cdot \leq y) \in\R_2$. If $b^{\pOrder + 1}(n b^d)^{(d+4)/(2d+4)}(\log n)^{-1/2}= O(1)$ and $(nb^d)^{-1}\log n = o(1)$, then
	\begin{align*}
		\sup_{\bw\in\W,y\in\reals} \big|\sqrt{n b^d}\big(\widehat{F}(y|\bw) - F(y|\bw)\big) - Z_n^R(\bw,y)\big|
		= O(\mathsf{r}_n) \quad \text{a.s.},
	\end{align*}
	where $\cov(Z_n^R(\bw_1,u_1),Z_n^R(\bw_2,u_2)) = n b^d \cov(\be_1^{\top}\bH_{\bw_1}^{-1}\mathbf{S}_{\bw_1,r_{u_1}}, \be_1^{\top} \bH_{\bw_2}^{-1}\mathbf{S}_{\bw_2,r_{u_2}})$ for all $(\bw_1,u_1),(\bw_2,u_2)$ in $\W \times \reals$ and $r_{u_1},r_{u_2}\in\R_2$.
\end{itemize}

This example gives a statistical application where Theorem \ref{thm: R-Process -- Main Theorem} offers a strict improvement on the accuracy of the Gaussian strong approximation over \cite{Rio_1994_PTRF}, and the improved Theorem \ref{thm: X-Process -- Main Theorem} upon incorporating a Lipschitz condition on the function class. See \cite[Section SA-IV.6]{Cattaneo-Yu_2024_AOS--SA} for omitted details. It remains an open question whether the result in this section provides the best Gaussian strong approximation for local polynomial regression or, more generally, for a local empirical process. The results presented are the best in the literature, but we are unaware of lower bounds that would confirm the approximation rates are unimprovable.

\section{Quasi-Uniform Haar Functions}\label{sec: Quasi-Uniform Haar Basis}

Assuming the existence of a surrogate measure and a normalizing transformation, or otherwise restricting the data generating process, Theorem \ref{thm: X-Process -- Main Theorem} established that the general empirical process \eqref{eq:X-process} indexed by VC-type Lipschitz functions can admit a strong approximation \eqref{eq:SA} at the optimal univariate KMT rate $\varrho_n = n^{-1/2} \log n$ when $d\in\{1,2\}$, and at the improved (but possibly suboptimal) rate $\varrho_n = n^{-1/d} \sqrt{\log n}$ when $d\geq3$, putting aside $\ttc_1$, $\ttc_2$, $\ttc_3$, $\ttM_{\scrH}$, $\ttL_{\scrH}$, $\ttTV_{\scrH}$, and $\mathtt{K}_{\scrH}$. The possibly suboptimal strong approximation rate arises from the $L_2$-approximation of the functions $h\in\scrH$ by a Haar basis expansion based on a carefully chosen \textit{dyadic} partition of a cover of $\X$. Likewise, Theorem \ref{thm: R-Process -- Main Theorem} established an improved uniform Gaussian strong approximation for the residual-based empirical process \eqref{eq:R-process}, but the result is also limited by the mean square projection error incurred by employing a Haar basis expansion based on a carefully chosen, asymmetric partitioning of the support of $\bz_i=(\bx_i,y_i)$.

Motivated by the limitations introduced by the mean square projection error underlying the proofs of Theorems \ref{thm: X-Process -- Main Theorem} and \ref{thm: R-Process -- Main Theorem}, this section presents uniform Gaussian strong approximations for $(X_n(h):h \in \scrH)$ and $(R_n(g,r):(g,r)\in\G\times\R)$ when $\scrH$ and $\G$ belong to the span of a Haar basis based on a \textit{quasi-uniform} partition with cardinality $L$, which can be viewed as an approximation based on $L\to\infty$ as $n\to\infty$. We do not require the existence of a normalizing transformation, allow for more general partitioning schemes than dyadic cells expansions, and impose minimal restrictions on the data generating process, while achieving the univariate KMT optimal strong approximation rate based on the effective sample size $n/L$ for all $d\geq1$. The strong approximation results presented in this section generalize two ideas from the regression Splines literature \citep{Huang_2003_AOS}: (i) the cells forming the Haar basis are assumed to be quasi-uniform with respect to a surrogate measure $\Q_\scrH$; and (ii) the number of active cells of the Haar basis affects the strong approximation. We apply the strong approximation results to histogram density estimation, and partitioning-based regression estimation based on Haar basis, which includes certain regression trees \citep{breiman1984} and other related methods \citep{Cattaneo-Farrell-Feng_2020_AOS}. Proof and omitted technical details are given in \cite[Section SA-V]{Cattaneo-Yu_2024_AOS--SA}.

\subsection{General Empirical Process}\label{sec: Quasi-Uniform Haar Basis -- X-process}

The following result is the analogue of Theorem \ref{thm: X-Process -- Main Theorem}.

\begin{thm}\label{thm: X-Process -- Approximate Dyadic Thm}
	Suppose $(\bx_i: 1 \leq i \leq n)$ are i.i.d. random vectors taking values in $(\reals^d, \mathcal{B}(\reals^d))$ with common law $\P_X$ supported on $\X \subseteq \reals^d$, and the following condition holds.
	\begin{enumerate}[label=(\roman*)]
		\item $\scrH \subseteq \operatorname{Span}\{\Indicator_{\Delta_l}: 0 \leq l < L\}$ is a class of Haar functions on $(\reals^d, \mathcal{B}(\reals^d),\P_X)$.
		\item There exists a surrogate measure $\bbQ_\scrH$ for $\P_X$ with respect to $\scrH$ such that $\{\Delta_l: 0 \leq l < L\}$ forms a \textit{quasi-uniform partition} of $\Q_\scrH$ with respect to $\bbQ_\scrH$:
		\begin{align*}
			\Q_\scrH \subseteq \sqcup_{0\leq l < L} \Delta_l \qquad\text{and}\qquad
			\frac{\max_{0 \leq l < L}\bbQ_\scrH(\Delta_l)}{\min_{0 \leq l < L}\bbQ_\scrH(\Delta_l)} \leq \rho < \infty.
		\end{align*}
		\item $\ttM_{\scrH} < \infty$.
	\end{enumerate}
	Then, on a possibly enlarged probability space, there exists a sequence of mean-zero Gaussian processes $(Z^X_n(h):h\in\scrH)$ with almost sure continuous trajectories on $(\scrH,\d_{\P_X})$ such that:
	\begin{itemize}
		\item $\E[X_n(h_1) X_n(h_2)] = \E[Z^X_n(h_1) Z^X_n(h_2)]$ for all $h_1, h_2 \in \scrH$, and
		\item $\P\big[\norm{X_n - Z^X_n}_{\scrH} > C_1 C_{\rho} \mathsf{P}_n(t)\big] \leq C_2 e^{-t} + L e^{-C_{\rho} n/L}$ for all $t > 0$,
	\end{itemize}
	where $C_1$ and $C_2$ are universal constants, $C_{\rho}$ is a constant that only depends on $\rho$, and
	\[\mathsf{P}_n(t) 
	= \min_{\delta\in(0,1)}\Big\{ \mathsf{H}_n(t,\delta) + \mathsf{F}_n(t,\delta) \Big\},\]
	with 
	\begin{align*}
		\mathsf{H}_n(t,\delta)
		& = \sqrt{\frac{\ttM_{\scrH} \ttE_{\scrH}}{n/L}} \sqrt{t + \log \ttN_{\scrH}(\delta,M_{\scrH})} \\
		&\qquad + \sqrt{\frac{\min\{\log_2 L,\mathtt{S}_{\scrH}^2\}}{n}} \ttM_{\scrH} (t + \log \ttN_{\scrH}(\delta,M_{\scrH})),
	\end{align*}
	where $\ttS_{\scrH} = \sup_{h\in\scrH} \sum_{l=1}^L \Indicator(\supp(h)\cap\Delta_l \neq \emptyset)$.
\end{thm}

This theorem shows that if $n^{-1}L\log(nL) \to 0$, then a valid strong approximation can be achieved with exponential probability concentration. The proof of Theorem \ref{thm: X-Process -- Approximate Dyadic Thm} leverages the fact that the $L_2$-projection error is zero by construction, but recognizes that \cite[Theorem 2.1]{Rio_1994_PTRF} does not apply because the partitions are \textit{quasi-dyadic}, preventing the use of the celebrated Tusn\'ady's inequality. Instead, in \cite{Cattaneo-Yu_2024_AOS--SA}, we present two technical results to circumvent that limitation: (i) we combine \cite[Lemma 2]{brown2010nonparametric} and \cite[Lemma 2]{sakhanenko1996estimates} to establish a version of Tusn\'ady's inequality that allows for more general binomial random variables $\mathsf{Bin}(n,p)$ with $\underline{p}\leq p \leq \overline{p}$, the error bound holding uniformly in $p$, as required by the quasi-dyadic partitioning structure; and (ii) we generalize \cite[Theorem 2.1]{Rio_1994_PTRF} to the case of quasi-dyadic cells.

Assuming a VC-type condition on $\scrH$, and putting aside $\ttM_{\scrH}$, $\ttE_{\scrH}$, and $\ttS_{\scrH}$, it follows that \eqref{eq:SA} holds with $\varrho_n = \sqrt{\log(n)} /\sqrt{n/L} + \log(n) /\sqrt{n}$. More generally, we have the following.

\begin{coro}[VC-type Haar Functions]\label{coro: X-process -- haar}
	Suppose the conditions of Theorem \ref{thm: X-Process -- Approximate Dyadic Thm} hold. In addition, assume that $\scrH$ is a VC-type class with function $\ttM_{\scrH}$ over $\Q_\scrH$ with constants $\ttc_{\scrH} \geq e$ and $\mathtt{d}_{\scrH} \geq 1$. Then, if $n^{-1}L\log(nL) \to 0$, \eqref{eq:SA} holds with
	\begin{equation*}
		\varrho_n = \sqrt{\frac{\ttM_{\scrH}\ttE_{\scrH}}{n/L}} \sqrt{\log n} + \sqrt{\frac{\min\{\log_2 L,\mathtt{S}_{\scrH}^2\}}{n}} \mathtt{M}_{\scrH} \log n.
	\end{equation*}
\end{coro}

We offer a simple statistical application of Theorem \ref{thm: X-Process -- Approximate Dyadic Thm} in the next example. 

\begin{example}[Histogram Density Estimation]\label{example: Histogram Density Estimation}
	The histogram density estimator of $f_X$ is
	\begin{align*}
		\check{f}_X(\bw) 
		= \frac{1}{n} \sum_{i = 1}^{n} \sum_{l = 0}^{P-1} \Indicator(\bw \in \Delta_l) \Indicator(\bx_i \in \Delta_l),
	\end{align*}
	where $\{\Delta_l: 0 \leq l < P\}$ are disjoint and satisfy $\max_{0 \leq l < P}\P_X(\Delta_l) \leq \rho \min_{0 \leq l < P}\P_X(\Delta_l).$
	
	For $L$ proportional to $\P_X(\Delta_l)^{-1}$, up to $\rho$, we establish a strong approximation for the localized empirical process $(\zeta_n(\bw):\bw\in\W)$, $\W\subseteq\X$, where
	\begin{align*}
		\zeta_n(\bw) = \sqrt{n L}\big(\check{f}_X(\bw) - \E[\check{f}_X(\bw)]\big) = X_n(h_{\bw}),
		\qquad  h_{\bw} \in \scrH,
	\end{align*}
	with $\scrH=\{h_{\bw}(\cdot) = L^{1/2} \sum_{l = 0}^{P-1} \Indicator(\bw \in \Delta_l) \Indicator(\cdot \in \Delta_l) : \bw\in\W \}$ a collection of Haar basis functions based on the partition $\{\Delta_l: 0 \leq l < P\}$. It follows that $\ttM_{\scrH,\reals^d} = L^{1/2}$ and $\ttS_{\scrH} = 1$.
	
	If $\W=\X$, then we set $L=P$, $\bbQ_\scrH=\P_X$, $\Q_\scrH=\X$, and the conditions of Theorem \ref{thm: X-Process -- Approximate Dyadic Thm} are satisfied with $\ttE_{\scrH} = L^{-1/2}$. Then, for $X_n=\zeta_n$, \eqref{eq:SA} holds with $\varrho_n = \log(nL)/\sqrt{n/L}$, assuming that $n^{-1}L\log(nL) \to 0$.

	If $\W\subsetneq \X$, assume $\W \subseteq \sqcup_{0 \leq l < P}\Delta_l$. If $\P_X(\sqcup_{0 \leq l < P} \Delta_l) < 1$, then $\{\Delta_l: 0 \leq l < P\}$ is no longer a quasi-uniform partition of $\X$ with respect to $\P_X$. The surrogate measure can help in this setting: we may add or refine cells to handle the residual probability $\P_X[(\sqcup_{0 \leq l < P} \Delta_l)^c]$. For example, suppose that for some $\mathring{P} \in \mathbb{N}$ we have 
	\begin{align*}
		\mathring{P} \leq \frac{\P_X((\sqcup_{0 \leq l < P} \Delta_l)^c)}{\min_{0 \leq l < P}\P_X(\Delta_l)}< \mathring{P} + 1.
	\end{align*}
	Set $L = P + \mathring{P}$. For any collection of disjoint cells $\{\Delta_l: P \leq l < L\}$ in $\X \cup \supp(\scrH)^c$, take $\bbQ_\scrH$ to agree with $\P_X$ on $\sqcup_{0 \leq l < P}\Delta_l$ and $\bbQ_\scrH(\Delta_l) = \mathring{P}^{-1}\P_X[(\sqcup_{0 \leq l < P} \Delta_l)^c]$ for $l=P,\dots, L-1$. Then, the enlarged class of cells $\{\Delta_l: 0 \leq l < L + K\}$ and the probability measure $\bbQ_\scrH$ satisfy conditions (i) and (ii) in Theorem~\ref{thm: X-Process -- Approximate Dyadic Thm}. It follows that $\ttE_{\scrH} = L^{-1/2}$ and hence, for $X_n=\zeta_n$, \eqref{eq:SA} holds with $\varrho_n = \log(nL)/\sqrt{n/L}$, assuming that $n^{-1}L\log(nL) \to 0$. In particular, the quasi-uniformity condition of $\P_X$ is required on a cover of $\W$, instead of on a cover of $\X$, at the expense of possibly increasing the number of cells to account for the residual probability $\P_X[(\sqcup_{0 \leq l < P}\Delta_l)^c]$.
\end{example}

Theorem \ref{thm: X-Process -- Approximate Dyadic Thm}, and in particular Example \ref{example: Histogram Density Estimation}, showcases the existence of a class of stochastic processes for which a uniform Gaussian strong approximation can be established with optimal univariate KMT rate in terms of the effective sample size $n/L$ for all $d\geq1$. This result is achieved because there is no projection error ($\scrH$ is spanned by a Haar basis), and the coupling error is controlled via our generalized Tusn\'ady's inequality. See \cite{Cattaneo-Yu_2024_AOS--SA} for details.

\subsection{Residual-Based Empirical Process}

The next result is the analogue of Theorem \ref{thm: R-Process -- Main Theorem}.

\begin{thm}\label{thm: R-Process -- Approximate Dyadic Thm}
	Suppose $(\bz_i=(\bx_i, y_i): 1 \leq i \leq n)$ are i.i.d. random vectors taking values in $(\reals^{d+1}, \mathcal{B}(\reals^{d+1}))$ with common law $\P_Z$, where $\bx_i$ has distribution $\P_X$ supported on $\X\subseteq\reals^d$, $y_i$ has distribution $\P_Y$ supported on $\Y\subseteq\reals$, and the following conditions hold.
	
	\begin{enumerate}[label=(\roman*)]
		\item $\G \subseteq \operatorname{Span}\{\Indicator_{\Delta_l}: 0 \leq l < L\}$ is a class of Haar functions on $(\reals^d, \mathcal{B}(\reals^d),\P_X)$.
		\item There exists a surrogate measure $\bbQ_\G$ for $\P_X$ with respect to $\G$ such that $\{\Delta_l: 0 \leq l < L\}$ forms a \textit{quasi-uniform partition} of $\Q_\G$ with respect to $\bbQ_\G$:
		\begin{align*}
			\Q_\G \subseteq \sqcup_{0\leq l < L} \Delta_l \qquad\text{and}\qquad
			\frac{\max_{0 \leq l < L}\bbQ_\G(\Delta_l)}{\min_{0 \leq l < L}\bbQ_\G(\Delta_l)} \leq \rho < \infty.
		\end{align*}
		\item $\G$ is a VC-type class with envelope function $\ttM_{\G}$ over $\Q_\G$ with $\ttc_{\G} \geq e$ and $\mathtt{d}_{\G} \geq 1$.
		
		\item $\R$ is a real-valued pointwise measurable class of functions on $(\reals, \mathcal{B}(\reals),\P_Y)$.
		
		\item $\R$ is a VC-type class with envelope $M_{\R,\Y}$ over $\Y$ with $\ttc_{\R,\Y}\geq e$ and $\mathtt{d}_{\R,\Y}\geq 1$, where $M_{\R,\Y}(y) + \mathtt{pTV}_{\R,(-|y|,|y|)} \leq \ttv (1 + |y|^{\alpha})$ for all $y \in \Y$, for some $\ttv>0$, and for some $\alpha\geq0$. Furthermore, if $\alpha>0$, then $\sup_{\bx \in \X}\E[\exp(|y_i|)|\bx_i = \bx] \leq 2$.
		
		\item There exists a constant $\ttk$ such that $|\log_2 \mathtt{E}_{\G}| + |\log_2 \ttM_{\G}| + |\log_2 L| \leq \ttk \log_2 n$.
		
	\end{enumerate}
	Then, on a possibly enlarged probability space, there exists a sequence of mean-zero Gaussian processes $(Z_n^R(g,r): (g,r)\in \G\times \R)$ with almost sure continuous trajectories on $(\G \times \R, \d_{\P_Z})$ such that:
	\begin{itemize}
		\item $\E[R_n(g_1, r_1) R_n(g_2, r_2)] = \E[Z^R_n(g_1, r_1) Z^R_n(g_2, r_2)]$ for all $(g_1, r_1), (g_2, r_2) \in \G \times \R$, and
		\item $\P[\norm{R_n - Z_n^R}_{\G \times \R} > C_1 C_{\ttv,\alpha} (C_{\rho} \mathsf{U}_n(t)+\mathsf{V}_n(t))] \leq C_2 e^{-t} + L e^{-C_{\rho}n/L}$ for all $t > 0$,
	\end{itemize}
	where $C_1$ and $C_2$ are universal constants, $C_{\ttv,\alpha} = \ttv \max\{1 + (2 \alpha)^{\frac{\alpha}{2}}, 1 + (4 \alpha)^{\alpha}\}$, $C_{\rho}$ is a constant that only depends on $\rho$,
	\begin{align*}
		\mathsf{U}_n(t)
		&= \Bigg( \sqrt{\frac{d \ttM_{\G} \ttE_{\G}}{n/L}} + \frac{\ttM_{\G}}{\sqrt{n}}(\log n)^{\alpha}\Bigg)
		(t + \ttk \log_2 n + \mathtt{d}\log (\ttc n))^{\alpha + 1}
	\end{align*}
	with $\ttc = \ttc_{\G} \ttc_{\R,\Y}$, $\mathtt{d} = \mathtt{d}_{\G} + \mathtt{d}_{\R,\Y}$,
	and
	\begin{align*}
		\mathsf{V}_n(t) = \Indicator(|\R|>1)\sqrt{\ttM_{\G} \ttE_{\G}} \Big(\max_{0 \leq l < L} \infnorm{\Delta_l}\Big) \ttL_{\mathscr{V}_{\R}} \sqrt{t + \ttk \log_2 n + \mathtt{d}\log(\ttc n)},
	\end{align*}
	with $\mathscr{V}_{\R} = \{\theta(\cdot,r) : r \in \R\}$, and $\theta(\bx,r) = \E[r(y_i)|\bx_i = \bx]$.
\end{thm}

The first term, $\mathsf{U}_n(t)$, can be interpreted as a ``variance'' contribution based on the effective sample size $n/L$, up to $\polylog(n)$ terms, while the second term, $\mathsf{V}_n(t)$, can be interpreted as a ``bias'' term that arises from the projection error for the conditional mean function $\E[r(y_i)|\bx_i = \bx]$, which may not necessarily lie in the span of Haar basis. In the special case when $\R$ is a singleton, we can construct the cells based on the condition distribution of $r(y_i) - \E[r(y_i)|\bx_i]$, thereby making the conditional mean function (and hence the ``bias'' term) zero, but such a construction is not possible when uniformity over $\R$ is desired.

Theorem \ref{thm: R-Process -- Approximate Dyadic Thm} gives the following uniform Gaussian strong approximation result.

\begin{coro}[VC-type Haar Basis]\label{coro: R-process -- haar}
	Suppose the conditions of Theorem \ref{thm: R-Process -- Approximate Dyadic Thm} hold with constants $\ttc$ and $\ttd$. Then, if $n^{-1}L\log(nL) \to 0$, $\norm{R_n - Z_n^R}_{\G \times \R} = O(\varrho_n)$ a.s. with
	\begin{align*}
		\varrho_n = 
		\sqrt{\frac{\ttM_{\G} \ttE_{\G}}{n/L}}(\log n)^{\alpha + 1} +  \frac{\ttM_{\G}}{\sqrt{n}}(\log n)^{2\alpha + 1} + \Indicator(|\R|>1) \sqrt{\ttM_{\G} \ttE_{\G}} (\max_{0 \leq l < L} \infnorm{\Delta_l}) \sqrt{\log n} .
	\end{align*}
\end{coro} 

Setting aside $\ttM_{\G}$ and $\ttE_{\G}$, an approximation rate is $(\log n)^{2\alpha + 1}(n/L)^{-1/2} + \Indicator(|\R|>1) (\max_{0 \leq l < L} \infnorm{\Delta_l}) \sqrt{\log n}$, which can achieve the optimal univariate KMT strong approximation rate based on the effective sample size $n/L$, up to a $\polylog(n)$ term, when $\R$ is a singleton function class. See \cite[Section SA-V]{Cattaneo-Yu_2024_AOS--SA} for details.

The next section illustrates Theorem \ref{thm: R-Process -- Approximate Dyadic Thm} with an example studying nonparametric regression estimation based on a Haar basis approximation. 

\subsection{Example: Haar Partitioning-based Regression}

Suppose $(\bz_i = (\bx_i, y_i), 1 \leq i \leq n)$ are i.i.d. random vectors taking values in $(\X \times \reals, \mathcal{B}(\X \times \reals))$ with $\X \subseteq \reals^d$. As in Section \ref{section: Local Polynomial Regression}, consider the regression estimand \eqref{eq: regression estimand}, focusing again on the two examples $\R_1$ and $\R_2$. Instead of local polynomial regression, we study the Haar partitioning-based estimator:
\begin{align*}
	\check{\theta}(\bw,r) = \bp(\bw)^{\top}\widehat{\bgamma}(r), \qquad \widehat{\bgamma}(r) = \argmin_{\bgamma \in \reals^L} \sum_{i = 1}^n \big(r(y_i) - \bp(\bx_i)^{\top}\bgamma\big)^2,
\end{align*}
where $\bp(\bu) = (\Indicator(\bu \in \Delta_l): 0 \leq l < L)$, and $\bw\in\W\subseteq\X$. As in Example \ref{example: Histogram Density Estimation}, either $\W=\X$ or $\W\subsetneq\X$, but for simplicity we discuss only the former case, and hence we assume that $\{\Delta_l: 0 \leq l < L\}$ is a quasi-uniform partition of $\Q_\scrH=\X$ with respect to $\bbQ_\scrH=\P_X$.

The estimation error can again be decomposed into three terms:
\begin{align*}
	&\check{\theta}(\bw,r) - \theta(\bw,r)\\
	&\quad = \underbrace{ \bp(\bw)^{\top} \bQ^{-1}\bT_{r} }_\text{linearization}
	+ \underbrace{ \bp(\bw)^{\top} (\widehat{\bQ}^{-1} - \bQ^{-1}) \bT_r }_\text{non-linearity error}
	+ \underbrace{ \E[\check{\theta}(\bw,r)|\bx_1,\cdots,\bx_n] - \theta(\bw,r) }_\text{smoothing bias},
\end{align*}
where $\bQ = \E[\bp(\bx_i) \bp(\bx_i)^{\top}]$, $\widehat{\bQ} = \frac{1}{n}\sum_{i = 1}^n \bp(\bx_i) \bp(\bx_i)^{\top}$, and $\bT_r = \frac{1}{n}\sum_{i =1}^n \bp(\bx_i) (r(y_i)-\E[r(y_i)|\bx_i])$. In this example, the linearization term takes the form
\begin{align*}
	\sqrt{n/L} \bp(\bw)^{\top} \bQ^{-1} \mathbf{T}_r
	= \frac{1}{\sqrt{n }} \sum_{i = 1}^n k_{\bw}(\bx_i) (r(y_i)-\E[r(y_i)|\bx_i])
	= R_n(g,r), \quad g \in \G, r \in \R_l,
\end{align*}
for $l = 1,2$, where $\G = \{k_{\bw}(\cdot): \bw \in \W\}$ with $k_{\bw}(\bu) = L^{-1/2} \sum_{0 \leq l < L} \Indicator(\bw \in \Delta_l) \Indicator(\bu \in \Delta_l)/\P_X(\Delta_l)$ the equivalent kernel. Under standard regularity conditions including smoothness and moment assumptions \citep[Section SA-V.3]{Cattaneo-Yu_2024_AOS--SA},
\begin{align*}
	&\sup_{r\in\R_1} \big|\be_1^{\top} (\widehat{\bQ}^{-1} - \bQ^{-1}) \bT_{r}\big|
	= O(\log(n L)L/n + (\log(n L)L/n)^{3/2} \log n) \qquad \text{a.s.},\\
	&\sup_{r\in\R_2} \big|\be_1^{\top} (\widehat{\bQ}^{-1} - \bQ^{-1}) \bT_{r}\big|
	= O(\log(n L)L/n) \qquad \text{a.s.},\\
	&\sup_{\bw\in\W,r\in\R_l} \big| \E[\check{\theta}(\bw,r)|\bx_1,\cdots,\bx_n] - \theta(\bw,r)\big|
	= O\big(\max_{0 \leq l <L} \infnorm{\Delta_l}\big) \qquad \text{a.s.}, \quad l=1,2,
\end{align*}
provided that $\log(nL)L/n \to0$. Finally, for the residual-based empirical process $(R_n(g,r): g \in \G, r \in \R_l)$, $l=1,2$, we apply Theorem \ref{thm: R-Process -- Approximate Dyadic Thm}. First, $\ttM_{\G} = L^{1/2}$ and $\ttE_{\G} = L^{-1/2}$, and we can take $\ttc_{\G} = L$ and $\mathtt{d}_{\G} = 1$ because $\G$ has finite cardinality $L$. For the singleton case $\R_1$, we can take $\ttc_{\R_1} = 1$ and $\mathtt{d}_{\R_1} = 1$, $\alpha = 1$ if $\sup_{\bx \in \X}\E[\exp(|y_i|)|\bx_i = \bx] \leq 2$, and condition (v) in Theorem~\ref{thm: R-Process -- Approximate Dyadic Thm} holds, which implies that $\norm{R_n - Z_n^R}_{\G \times \R_1} = O(\varrho_n)$ a.s. with
\begin{equation*}
	\varrho_n =  \frac{\log(n L)^2}{\sqrt{n/L}},
\end{equation*}
provided that $\log(nL)L/n \to0$. For the VC-Type class $\R_2$, we can verify condition (v) in Theorem~\ref{thm: R-Process -- Approximate Dyadic Thm} with $\alpha=0$, and we can take $\ttc_{\R_2}$ to be some universal constant and $\mathtt{d}_{\R_2} = 2$ by \cite[][Theorem 2.6.7]{wellner2013weak}, which implies that $\norm{R_n - Z_n^R}_{\G \times \R_1} = O(\varrho_n)$ a.s. with
\begin{equation*}
	\varrho_n = \frac{\log (n L)}{\sqrt{n/L}} + \max_{0 \leq l < L} \infnorm{\Delta_l},
\end{equation*}
provided that $\log(n)L/n \to0$. A uniform Gaussian strong approximation for the Haar partitioning-based regression processes $(\sqrt{n/L}(\check{\theta}(\bw,r) - \theta(\bw,r)) : (\bw,r) \in \W \times \R_l)$, $l = 1,2$, follows directly from the results obtained above, as illustrated in Section \ref{section: Local Polynomial Regression}.

This example showcases a statistical application of our strong approximation result (Theorem \ref{thm: R-Process -- Approximate Dyadic Thm}) where the optimal univariate KMT strong approximation rate based on the effective sample size $n/L$ is achievable, up to $\polylog(n)$ terms and the complexity of $\R$. See \cite[Section SA-V.3]{Cattaneo-Yu_2024_AOS--SA} for omitted details.

\begin{acks}[Acknowledgments]
    We specially thank Boris Hanin for many insightful discussions. We also thank Rajita Chandak, Jianqing Fan, Kengo Kato, Jason Klusowski, Xinwei Ma, Boris Shigida, Jennifer Sun, Rocio Titiunik, Will Underwood, and two reviewers for their comments and suggestions.
\end{acks}
\begin{funding}
    The first author was supported by the National Science Foundation through grants DMS-2210561 and SES-2241575.
\end{funding}

\begin{supplement}
  \stitle{Proofs and other technical results}
  \sdescription{The supplementary material
    \citep{Cattaneo-Yu_2024_AOS--SA}
    collects detailed proofs
    of our main results, and also provides other technical results that may
  be of independent interest.}
\end{supplement}


\bibliographystyle{imsart-nameyear}
\bibliography{CY_2024_AOS--bib.bib}       


\end{document}